\documentclass[a4paper,11pt]{article}
\usepackage{amsmath,amssymb,amsthm,a4wide,mathrsfs}
\usepackage[utf8]{inputenc}
\usepackage[T1]{fontenc}
\usepackage{stmaryrd}
\usepackage{color,graphicx,subcaption}

\newcommand{\nc}{\newcommand}
\nc{\dsp}{\displaystyle}
\nc{\RR}{\mathbb{R}}
\nc{\CC}{\mathbb{C}}
\nc{\mrm}{\mathrm}
\nc{\bx}{\boldsymbol{x}}
\nc{\by}{\boldsymbol{y}}
\nc{\bn}{\boldsymbol{n}}
\nc{\loc}{\mrm{loc}}
\nc{\mH}{\mrm{H}}
\nc{\mL}{\mrm{L}}
\nc{\mB}{\mrm{B}}
\nc{\mA}{\mrm{A}}
\nc{\mJ}{\mrm{J}}
\nc{\mT}{\mrm{T}}
\nc{\mZ}{\mrm{Z}}
\nc{\mR}{\mrm{R}}
\nc{\mQ}{\mrm{Q}}
\nc{\mS}{\mrm{S}}
\nc{\mD}{\mrm{D}}
\nc{\mP}{\mrm{P}}
\nc{\mV}{\mrm{V}}
\nc{\mK}{\mrm{K}}
\nc{\mW}{\mrm{W}}

\nc{\mbH}{\mathbb{H}}
\nc{\mbX}{\mathbb{X}}
\nc{\dir}{\textsc{d}}
\nc{\neu}{\textsc{n}}
\nc{\lbr}{\lbrack}
\nc{\rbr}{\rbrack}
\nc{\ctru}{\mathfrak{u}}
\nc{\ctrv}{\mathfrak{v}}
\nc{\ctrw}{\mathfrak{w}}
\nc{\ctrp}{\mathfrak{p}}
\nc{\ctrq}{\mathfrak{q}}
\nc{\ctrf}{\mathfrak{f}}
\nc{\ctrg}{\mathfrak{g}}
\nc{\ctrh}{\mathfrak{h}}

\nc{\calC}{\mathcal{C}}

\nc{\Id}{\mrm{Id}}
\nc{\llbr}{\llbracket}
\nc{\rrbr}{\rrbracket}
\nc{\vphi}{\varphi}
\nc{\bvphi}{\boldsymbol{\varphi}}
\nc{\bp}{\boldsymbol{p}}
\nc{\green}{\mathscr{G}}
\nc{\SL}{\mrm{SL}}
\nc{\DL}{\mrm{DL}}
\nc{\llangle}{\langle\!\langle}
\nc{\rrangle}{\rangle\!\rangle}
\nc{\bpsi}{\boldsymbol{\psi}}

\renewcommand{\div}{\mrm{div}}

\newtheorem{proposition}{Proposition}[section]
\newtheorem{theorem}    {Theorem}    [section]
\newtheorem{lemma}      {Lemma}      [section]
\newtheorem{corollary}  {Corollary}  [section]

\title{\textbf{A new variant of the Optimised Schwarz Method for} \\
  \textbf{arbitrary non-overlapping subdomain partitions}}\date{}
\author{X.Claeys\footnote{Sorbonne Université, Université Paris-Diderot SPC, CNRS, Inria,
  Laboratoire Jacques-Louis Lions, équipe Alpines, email: claeys@ann.jussieu.fr}}

\begin{document}

\maketitle

\begin{abstract}
  We consider a scalar wave propagation in harmonic regime
  modelled by Helmholtz equation with heterogeneous coefficients.  
  Using the Multi-Trace Formalism (MTF), we propose a
  new variant of the Optimized Schwarz Method (OSM) that
  can accomodate the presence of cross-points in the subdomain
  partition. This leads to the derivation of a strongly coercive
  formulation of our Helmholtz problem posed on the union of all interfaces.
  The corresponding operator takes the form "identity + contraction".
\end{abstract}

\section{Introduction}
The effective solution to large scale wave propagation problems
relates to a wide range of applications and yet remains a challenge,
in particular when simulating highly oscillatory phenomena.
With the growing importance of parallel computing, an intense
research effort has been dedicated, in recent years, to the
development of domain decomposition strategies that can be efficiently
applied to wave propagation problems.

There is now a vast litterature and a rich arsenal of well established
domain decomposition techniques to deal with symmetric positive
problems see e.g. \cite{zbMATH02113718,zbMATH06029154,zbMATH06534518}.
By essence though, wave propagation does not fall into this symmetric positive
framework and domain decomposition is much less developped for waves,
from the point of view of both theory and effective numerical computation.
\quad\\

\begin{figure}[h]
  \hspace{0.25cm}
  \begin{subfigure}{.3\textwidth}
    \includegraphics[height=3cm]{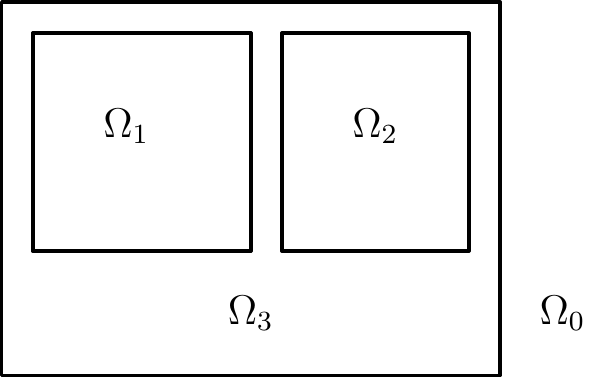}
    \caption{}\label{subfig1}
  \end{subfigure}
  \hspace{0.75cm}
  \begin{subfigure}{.3\textwidth}
    \includegraphics[height=3cm]{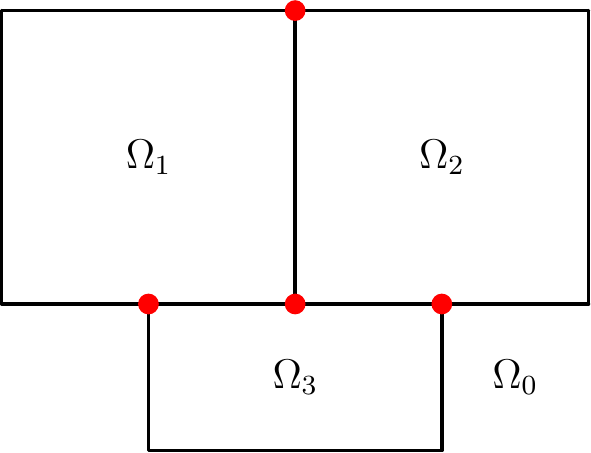}
    \caption{}\label{subfig2}
  \end{subfigure}
  \begin{subfigure}{.3\textwidth}
    \includegraphics[height=3cm]{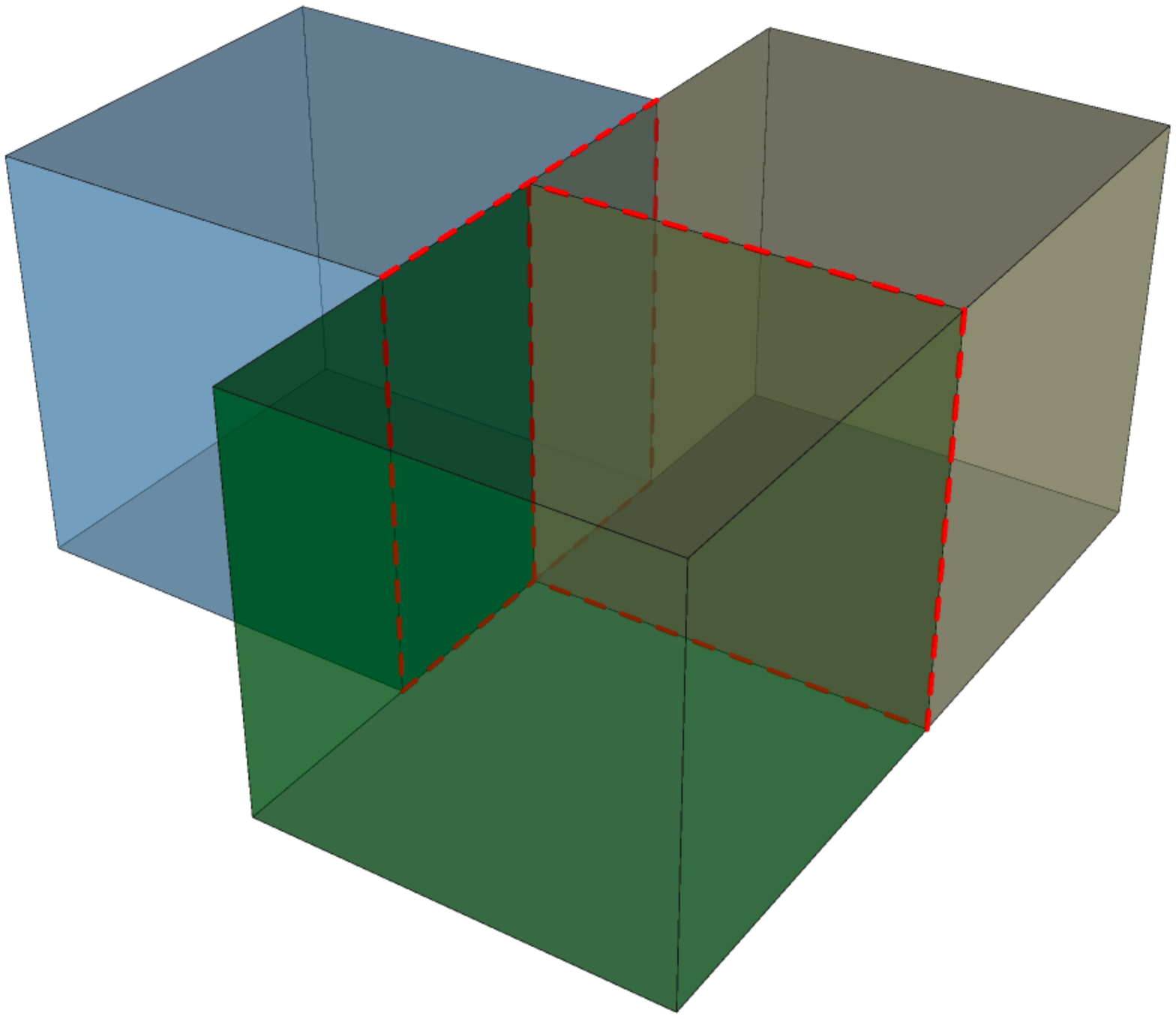}
    \caption{}\label{subfig3}
  \end{subfigure}
  \caption{Examples of subdomain partitions in 2D (a \& b) and 3D (c)
    with 4 subdomains (3 bounded + exterior). There is no cross point in (a),
    and cross points are red dots in (b) and red dashed lines in (c).}\label{Fig1}
\end{figure}

\noindent 
In the case of harmonic regime propagation, the Optimized Schwarz Method (OSM)
appears to be one of the most effective available approaches for domain decomposition
in a wave context. A general overview of this method and its numerous variants is given in
\cite{zbMATH07020343}.  In OSM, the coupling of subdomains is maintained through
transmission conditions at interfaces, and these transmission conditions are formulated
in terms of ingoing and outgoing trace operators involving impedance coefficients.
The efficiency of OSM crucially depends on the choice of these impedances.

The Optimized Schwarz Method was originally introduced in \cite{MR1071633,MR1105979,MR1291197,MR1227838}
considering general non-overlapping partition of the computational domain and constant scalar
impedance coefficients. Although, in such a general geometrical setting, OSM with scalar impedance
was proved to converge, no assessment was provided as regards the rate of convergence. In practice,
the convergence could be slow. This was improved by Collino and Joly in
\cite{MR1764190,LECOUVEZ2014403, 10.1007/978-3-319-93873-8_16} where the authors proposed
operator valued self-adjoint positive impedance coefficients and could establish geometric convergence
of the method assuming that the subdomain partition does not involve any cross point i.e. point of adjacency
of three interfaces (or one interface meeting the boundary of the compuational domain),
see Fig.\ref{Fig1} above. In another series of contributions Antoine, Geuzaine
and their collaborators  \cite{zbMATH06732120,zbMATH06660584,zbMATH06666466,zbMATH06039354,zbMATH07007340}
considered the case of impedance coefficients approaching appropriate
Dirichlet-to-Neumann maps and obtained fastly converging numerical methods. Here also,
the numerical methods were observed to be of good quality only when the subdomain partition
does not contain any cross-point. 

While much litterature has then been dedicated to the question of how to choose impedance
coefficients, cross points remained a thorny issue which, recently, has received a renewed
attention \cite{modave:hal-01925160}.  A very similar issue related to cross-points also
arises in a different context: the derivation of Boundary Integral Equations (BIE) adapted to
multi-domain scattering. The Multi-Trace Formalism (MTF) was introduced in \cite{MR3069956,MR3403719,
  zbMATH06286985,zbMATH06185811} as a complete
framework for dealing with multi-domain BIE. From the perspective of functional analysis, MTF offers
a clean treatment of cross-points. It would thus appear natural to try using the techniques developped
in the Multi-Trace framework for dealing properly with cross points in Optimized Schwarz domain
decomposition. This is precisely the aim of the present contribution.

\quad\\
In the present article, we introduce a new variant of the Optimized Schwarz Method for the solution of
Helmholtz equation with heterogeneous material coefficients through Formulation (\ref{OSMFormulation}).
This new variant can be applied with any non-overlapping partition of the propagtion
medium into Lipschitz subdomains, no matter the presence of cross-points. The operator of the
corresponding formulation takes the form "$\text{identity} + \text{contraction}$" in an appropriate trace space,
and we show that this operator is coercive. The key ingredient in this formulation is a non-local
exchange operator used to enforce transmission conditions. Such exchange operator has always
existed in previous versions of OSM, but it was so far systematically assumed to be a local
operator consisting in swapping the traces from both sides of each interface of the subdomain
partition. The exchange operator we consider here is more elaborate, which is the main novelty
of our approach.

It should be mentionned that the present contribution is purely analytical and that, in its present
form, this new variant of OSM does not seem appropriate for actual numerical computations. This is
why we do not report on numerical results. In a forthcoming article we will propose a discrete version
of the present formulation that is better suited for numerics. We still believe that the formulation
we present here is an interesting theoretical object. In particular, it yields a strongly coercive
formulation of Helmholtz problem which is not trivial: the derivation of coercive formulations for
Helmoltz equation has been, in itself, the subject of recent attention \cite{zbMATH06308881}.
In addition, in the case of piecewise constant material cooefficients, Formulation
(\ref{OSMFormulation}) can also be used as a multi-domain coupling scheme for the solution to
scattering problems by means of boundary integral formulation. In the particular case of piecewise constant
coefficients, the new formulation presented here can be considered as an alternative to
other multi-domain BIE such as Multi-Trace \cite{zbMATH06286985}, Boundary Element Tearing and Interconnecting
\cite{zbMATH02038442}, or Rumsey's reaction principle \cite{zbMATH04136456}.

\section{Geometry and problem under study}
In the present article, we are interested in a classical wave propagation
problem in harmonic regime set in an heterogeneous medium in $\RR^{d}$ for $d=1,2$ or $3$. 
We consider two essentially bounded measurable functions $\mu:\RR^{d}\to \RR_{+}$ and
$\kappa:\RR^{d}\to \CC_{+}$, and we assume that there exist constants $\kappa_{0},\rho_{0}>0$
such that 
\begin{equation}\label{HypoMat1}
  \begin{array}{rl}
    \textit{i)}   & \sup_{\bx\in \RR^{d}}(\vert\mu(\bx)\vert+\vert \mu^{-1}(\bx)\vert +
                    \vert \kappa(\bx)\vert ) <+\infty\\[5pt]
    \textit{ii)}  & \Re e\{\kappa(\bx)\}\geq 0,\; \Im m\{\kappa(\bx)\}\geq 0,\;\kappa(\bx)\neq 0
                    \quad \forall \bx\in \RR^{d}\\[5pt]
    \textit{iii)} & \kappa(\bx) = \kappa_{0}\;\text{and}\;\mu(\bx) = 1
                    \;\; \text{for}\;\vert\bx\vert>\rho_{0}
  \end{array}
\end{equation}
These assumptions are rather general yet reasonable enough to make the scattering problem
we wish to examine properly well posed. We insist that we do not assume $\kappa,\mu$ to be
piecewise constant. For some continuous functional $f\in \mL^{2}(\RR^{d})$ with bounded support,
we wish to solve the problem
\begin{equation}\label{IntialPb}
  \left\{\begin{aligned}
    & u\in\mH^{1}_{\loc}(\RR^{d})\;\text{such that}\\[3pt]
    & -\div(\mu\nabla u) - \kappa^{2} u = f\quad \text{in}\;\RR^{d},\\
    & \lim_{\rho\to\infty}\int_{\partial\mB_{\rho}}\vert \partial_{\rho}u-
    \imath\kappa_{0}u\vert^{2}d\sigma_{\rho} = 0.
  \end{aligned}\right.
\end{equation}
where $\mB_{\rho}$ refers to the ball of radius $\rho$ centered at $0$, $\sigma_{\rho}$
is the associated surface measure, and $\partial_{\rho}$ is the partial derivative
with respect to $\vert\bx\vert$. Well-posedness of the problem above is a classical
result of scattering theory, see e.g. \cite[Chap.3]{MR841971} or \cite[Chap.7]{MR2986407}.

\quad\\
We wish to solve this problem by means of non-overlapping Domain
Decomposition (DDM), which leads us to introduce a subdomain partitionning $\RR^{d} = \cup_{j=0}^{\mJ}
\overline{\Omega}_{j}$ with $\Omega_{j}\cap \Omega_{k} = \emptyset$ if $j\neq k$, each $\Omega_{j}$ is a
Lipschitz domain, and $\Omega_{j}$ is bounded for $j\neq 0$. The "skeleton" will refer to the union
of all interfaces between subdomains
\begin{equation*}
\Gamma = \partial\Omega_{0}\cup \dots\cup \partial\Omega_{\mJ}.
\end{equation*}
We emphasize that such geometrical configuration allows the presence
of junction points i.e. points where three subdomains or more abut. Examples of such
non-overlapping multi-domain configurations are given in Fig.\ref{Fig1}.

For the sake of simplicity, we make further regularity assumptions on material coefficients
in each subdomain, assuming that $\mu$ is Lipschitz regular in each subdomain,
\begin{equation}\label{HypoMat2}
  \begin{aligned}
    & \nabla\mu_{j}\in \mL^{\infty}(\Omega_{j})\;\forall j=0\dots \mJ,\\
    & \text{where}\quad \mu_{j}:=\mu\vert_{\Omega_{j}}.
  \end{aligned}
\end{equation}
Assumptions (\ref{HypoMat1})-(\ref{HypoMat2}) allow the coefficients $\mu,\kappa$
to jump across the interfaces $\partial\Omega_{j}\cap\partial\Omega_{k}$, but discards jumps of $\mu$ inside each
subdomain. In particular, this setting includes the case where $\mu,\kappa$
are piecewise constant with respect to the subdomain partition. 

\quad\\
Problem (\ref{IntialPb}) can be decomposed according to the subdomain
partition introduced above, leading to wave equations in each subdomain coupled by
transmission conditions imposed through each interface 
\begin{align}
  & \left\{\begin{aligned}
      & u\in\mH^{1}_{\loc}(\overline{\Omega}_{j})\;\text{such that}\\[3pt]
      & -\div(\mu\nabla u) - \kappa^{2} u = f\quad \text{in}\;\Omega_{j},\\
      & \lim_{\rho\to\infty}\int_{\partial\mB_{\rho}}\vert \partial_{\rho}u-
      \imath\kappa_{0}u\vert^{2}d\sigma_{\rho} = 0,
    \end{aligned}\right. \label{IntialP2:1}\\[10pt]
  & \left\{\begin{aligned}
      & u\vert_{\partial\Omega_{j}}^{\mrm{int}} - u\vert_{\partial\Omega_{k}}^{\mrm{int}} = 0
      & \forall j,k = 0\dots n\\[3pt]
      & \mu_{j}\partial_{n_{j}}u\vert_{\partial\Omega_{j}}^{\mrm{int}} + \mu_{k}\partial_{n_{k}}u\vert_{\partial\Omega_{k}}^{\mrm{int}} = 0
      & \text{on}\;\partial\Omega_{j}\cap \partial\Omega_{k}.\\
    \end{aligned}\right. \label{IntialP2:2}
\end{align}
where $\bn_{j}$ refers to the normal vector field on $\partial\Omega_{j}$ directed toward
the exterior of $\Omega_{j}$, and $\partial_{n_{j}}v:=\bn_{j}\cdot\nabla v$. The boundary
traces coming into play in the transmission conditions above are taken from the interior
of the subdomains, which is the meaning of the "\textrm{int}" superscript.

\quad\\
The present contribution will consist in deriving a strongly coercive reformulation
of Problem (\ref{IntialP2:1})-(\ref{IntialP2:2}) of the form "identity+contraction". This
reformulation will be posed in a space of trace on the skeleton $\Gamma$.

\section{Trace spaces and operators}
The treatment of interfaces between subdomains is a crucial aspect of any
domain decomposition strategy, both for constructing or analysing it. As a
consequence we pay a special attention to trace spaces.

\subsection{Volume based spaces}\label{VolumeSpaces}
First of all we need to fix a few notations related to classical
volume based function spaces. For any Lipschitz domain $\Omega\subset \RR^{d}$,
the space $\mL^{2}(\Omega)$ will refer to square integrable functions equipped with the norm
$\Vert \vphi\Vert_{\mL^{2}(\Omega)}^{2} := \int_{\Omega}\vert \vphi\vert^{2}d\bx$.
The Sobolev space $\mH^{1}(\Omega):=\{\vphi\in\mL^{2}(\Omega),\;\nabla\vphi\in
\mL^{2}(\Omega)^{d}\;\}$ will be equipped with the norm
\begin{equation}\label{H1norm}
  \Vert v\Vert_{\mH^{1}(\Omega)}^{2} := \Vert \nabla v\Vert_{\mL^{2}(\Omega)}^{2} + \gamma^{-2}\Vert v\Vert_{\mL^{2}(\Omega)}^{2}
\end{equation}
In this definition $\gamma>0$ refers to a parameter that will be fixed all through this article.
Occasionally we shall consider $\mH(\div,\Omega):=\{\bpsi\in \mL^{2}(\Omega)^{d}, \div(\bpsi)
\in \mL^{2}(\Omega)\}$ and $\mH^{1}(\Delta,\Omega):=\{\vphi\in\mH^{1}(\Omega),\;\Delta\vphi
\in\mL^{2}(\Omega)\}$ equipped with the norm given by $\Vert \vphi\Vert_{\mH^{1}(\Delta,\Omega)}^{2}
:= \Vert \vphi\Vert_{\mH^{1}(\Omega)}^{2}+\Vert \Delta \vphi\Vert_{\mL^{2}(\Omega)}^{2}$.
Finally if $\mH(\Omega)$ refers to any of the spaces introduced above, then
$\mH_{\loc}(\overline{\Omega})$
shall refer to all functions $v:\Omega\to \CC$ such that $v\vphi\in \mH(\Omega)$ for all
$\vphi\in \mathscr{C}^{\infty}_{\mrm{comp}}(\RR^{d}):=\{ \psi\in\mathscr{C}^{\infty}(\RR^{d}), \;
\mrm{supp}(\psi)\;\text{bounded}\}$.

\subsection{Traces on the boundary of a single subdomain}
For any Lipschitz open set $\Omega\subset \RR^{d}$,
we shall refer to the space of Dirichlet traces $\mH^{1/2}(\partial\Omega):=
\{v\vert_{\partial\Omega}, v\in \mH^{1}(\Omega)\}$ equipped with the norm 
\begin{equation}\label{TraceDirNorm}
\Vert v\Vert_{\mH^{1/2}(\partial\Omega)}:=\min\{\Vert \varphi\Vert_{\mH^{1}(\Omega)}, \varphi\vert_{\partial\Omega} = v\}.
\end{equation}
The space of Neumann traces $\mH^{-1/2}(\partial\Omega)$
will be defined as the dual to $\mH^{1/2}(\partial\Omega)$ equipped with the corresponding
canonical dual norm $\Vert p\Vert_{\mH^{-1/2}(\partial\Omega)}:=\sup_{v\in \mH^{1/2}(\partial\Omega)}
\vert\langle p,v\rangle_{\partial\Omega}\vert/\Vert v\Vert_{\mH^{1/2}(\partial\Omega)}$. Here
$v\mapsto \langle p,v\rangle_{\partial\Omega}:= p(v)$ simply refers to the action of $p$ on $v$,
so that $(p,v)\mapsto \langle p,v\rangle_{\partial\Omega}$ is a bilinear (not sesquilinear) form.
As regards duality pairing, we shall also equivalently write $\langle v,p\rangle_{\partial\Omega} :=
\langle p,v\rangle_{\partial\Omega}$ and
\begin{equation*}
\int_{\partial\Omega}pv d\sigma=\langle p,v\rangle_{\partial\Omega}.
\end{equation*}
We will also equip the space of pairs of Dirichlet/Neumann traces
with its own duality pairing. Although many choices are possible, we use a
skew-symmetric pairing that appears naturally in energy conservation calculus,
defined by
\begin{equation}
  \begin{aligned}
    &\lbr (u,p), (v,q)\rbr_{\partial\Omega}:=\langle u,q\rangle_{\partial\Omega}
     -\langle v,p\rangle_{\partial\Omega}\\
    & \text{for}\; (u,p)\;\text{and}\; (v,q)
    \;\text{in}\;
    \mH^{+\frac{1}{2}}(\partial\Omega)\times \mH^{-\frac{1}{2}}(\partial\Omega).
  \end{aligned}
\end{equation}
Note that this pairing does not involve any complex conjugation.
Let $\bn_{\Omega}$ refer to the normal vector field on $\partial\Omega$ 
directed toward the exterior of $\Omega$.
Each Lipschitz open set $\Omega\subset \RR^{d}$ with bounded
boundary gives rise to continuous operators $\tau_{\dir}^{\Omega}:
\mH^{1}_{\loc}(\Omega)\to \mH^{1/2}(\partial\Omega)$,
$\tau^{\Omega}_{\neu}:\mH^{1}_{\loc}(\Delta,\overline{\Omega})\to \mH^{-1/2}(\partial\Omega)$ and
$\tau^{\Omega}:\mH^{1}_{\loc}(\Delta,\overline{\Omega})\to\mH^{1/2}(\partial\Omega)\times 
\mH^{-1/2}(\partial\Omega)$ uniquely defined by
\begin{equation}\label{InteriorTraceOperator}
  \begin{aligned}
    & \tau_{\dir}^{\Omega}(\varphi):= \varphi\vert_{\partial\Omega}
    \quad \text{and}\quad
    \tau_{\neu}^{\Omega}(\varphi):= \bn_{\Omega}\cdot\nabla \varphi\vert_{\partial\Omega},\\[5pt]
    & \tau^{\Omega}(\varphi):=(\tau^{\Omega}_{\dir}(\varphi),\tau^{\Omega}_{\neu}(\varphi)) 
    \quad\quad \forall\varphi\in\mathscr{C}^{\infty}(\overline{\Omega}).
  \end{aligned}
\end{equation}

\subsection{Scalar products and Dirichlet-to-Neumann maps}
For any $v\in \mH^{1/2}(\partial\Omega)$ let $\phi_{\dir}(v)\in \mH^{1}(\Omega)$
refer to the unique element that achieves the minimum in (\ref{TraceDirNorm}) i.e. such that
$\Vert v\Vert_{\mH^{1/2}(\partial\Omega)} =
\Vert \phi_{\dir}(v)\Vert_{\mH^{1}(\Omega)}$. Writing Euler's identity for this minimisation problem,
we see that $\int_{\Omega}\nabla \phi_{\dir}(v)\cdot\nabla \vphi +\gamma^{-2}\phi_{\dir}(v)\vphi\, d\bx = 0\;
\forall \vphi\in \mH^{1}_{0}(\Omega)$, which re-writes $-\Delta \phi_{\dir}(v)+\gamma^{-2}\phi_{\dir}(v) = 0$ in
$\Omega$. Then we introduce a so-called Dirichlet-to-Neumann (DtN) map $\mT_{\Omega}:=
\tau^{\Omega}_{\neu}\cdot\phi_{\dir}:\mH^{1/2}(\partial\Omega)\to
\mH^{-1/2}(\partial\Omega)$. To be more explicit $\mT_{\Omega}$ is defined by 
\begin{equation}\label{PDEDtN}
  \begin{aligned}
    & \mT_{\Omega}(v):=\bn_{\Omega}\cdot\nabla \phi_{\dir}(v)\vert_{\partial\Omega}\\
    & \text{where}\;\phi_{\dir}(v)\in \mH^{1}(\Omega)\;\text{satisfies}\\
    & \Delta \phi_{\dir}(v)-\gamma^{-2}\phi_{\dir}(v) = 0\quad \text{in}\;\Omega\\
    & \phi_{\dir}(v)\vert_{\partial\Omega} = v\quad\text{on}\;\partial\Omega.
  \end{aligned}
\end{equation}
This DtN map actually induces the scalar product associated to the norm (\ref{TraceDirNorm}).
First of all observe that $\overline{\phi_{\dir}(u)} = \phi_{\dir}(\overline{u})$ obviously. 
Next, according to the PDE satisfied by $\phi_{\dir}$ in (\ref{PDEDtN}), applying Green's
formula we obtain $\int_{\Omega}\nabla \phi_{\dir}(u)\cdot \nabla \phi_{\dir}(\overline{v}) +
\gamma^{-2}\phi_{\dir}(u) \phi_{\dir}(\overline{v}) d\bx =
\int_{\partial\Omega}\phi_{\dir}(\overline{v})\,\bn_{\Omega}\cdot\nabla \phi_{\dir}(u) d\sigma =
\langle \mT_{\Omega}(u),\overline{v}\rangle_{\partial\Omega}$. From this calculus it
is clear that  $\langle \mT_{\Omega}(u),\overline{v}\rangle_{\partial\Omega} = \langle
\mT_{\Omega}(\overline{v}),u\rangle_{\partial\Omega}$. Since by the defintion of $\phi_{\dir}$ we have
$\Vert u\Vert_{\mH^{1/2}(\partial\Omega)} = \Vert \phi_{\dir}(u)\Vert_{\mH^{1}(\Omega)}$, we
can take the following as scalar product on the Dirichlet trace spaces
\begin{equation}\label{ScalarProductDirTrace}
  (u,v)_{\mH^{1/2}(\partial\Omega)}:=\langle\mT_{\Omega}(u),\overline{v}\rangle_{\partial\Omega}
  \quad\quad\text{for}\;u,v\in\mH^{1/2}(\partial\Omega). 
\end{equation}
According to Riesz representation
theorem, for any $p\in\mH^{-1/2}(\partial\Omega)$ there exists a unique
$\vphi_{p}\in\mH^{1/2}(\partial\Omega)$ such that $\langle p,v\rangle_{\partial\Omega} =
(\vphi_{p},v)_{\mH^{1/2}(\partial\Omega)}= \langle \mT_{\Omega}(\vphi_{p}),v\rangle$ for all
$v\in\mH^{1/2}(\partial\Omega)$. Hence $\vphi_{p} = (\mT_{\Omega})^{-1}(p)$ and
$\Vert p\Vert_{\mH^{-1/2}(\partial\Omega)}^{2} =
\Vert \vphi_{p}\Vert_{\mH^{1/2}(\partial\Omega)}^{2} = \langle \mT_{\Omega}(\vphi_{p}),
\overline{\vphi}_{p} \rangle_{\partial\Omega} = \langle p,\mT_{\Omega}^{-1}(\overline{p}) \rangle_{\partial\Omega}$. As a
consequence the norm on Neumann data is induced by the following scalar product
\begin{equation}\label{NeumannScalarProduct}
  (p,q)_{\mH^{-1/2}(\partial\Omega)}:=\langle p,\mT_{\Omega}^{-1}(\overline{q})\rangle_{\partial\Omega}
    \quad\quad\text{for}\;p,q\in\mH^{-1/2}(\partial\Omega). 
\end{equation}

\subsection{Traces in a multi-domain setting}\label{TracesMultiDomain}
We will also need to consider cartesian products of Dirichlet or Neumann
trace spaces based on the boundary of each subdomain of the partition, which we call
multi-trace spaces defined as follows
\begin{equation}
  \begin{aligned}
    & \mbH_{\dir}(\Gamma):=\mH^{+\frac{1}{2}}(\partial\Omega_{0})
    \times\dots\times \mH^{+\frac{1}{2}}(\partial\Omega_{\mJ}),\\[-3pt]
    & \mbH_{\neu}(\Gamma):=\mH^{-\frac{1}{2}}(\partial\Omega_{0})\times\dots\times 
    \mH^{-\frac{1}{2}}(\partial\Omega_{\mJ}),\\[-3pt]
    & \mbH(\Gamma)_{\textcolor{white}{\dir}}:=\Pi_{j=0\dots \mJ}
    \mH^{+\frac{1}{2}}(\partial\Omega_{j})\times \mH^{-\frac{1}{2}}(\partial\Omega_{j}).
  \end{aligned}
\end{equation}
equipped with  $\Vert \ctrp\Vert_{\mbH_{\neu}(\Gamma)}^{2} := \Vert p^{0}\Vert_{\mH^{-1/2}(\partial\Omega_{0})}^{2}
+ \dots+ \Vert p^{\mJ}\Vert_{\mH^{-1/2}(\partial\Omega_{\mJ})}^{2}$ for
$\ctrp = (p^{j})_{j=0}^{\mJ}\in\mbH_{\neu}(\Gamma)$, and analogous definitions for
$\Vert \;\;\Vert_{\mbH_{\dir}(\Gamma)}$ and  $\Vert \;\;\Vert_{\mbH(\Gamma)}$.
The multi-trace space $\mbH(\Gamma)$ coincides with $\mbH_{\dir}(\Gamma)\times\mbH_{\neu}(\Gamma)$
through a re-ordering of traces which is why, when considering an element $\ctru =
(u^{j}_{\dir},u^{j}_{\neu})_{j=0}^{\mJ}\in \mbH(\Gamma)$, 
we will sometimes commit a slight abuse of notation writing "$\ctru = (\ctru_{\dir},\ctru_{\neu})$" to refer to
the Dirichlet components $\ctru_{\dir} = (u^{j}_{\dir})_{j=0}^{\mJ}\in\mbH_{\dir}(\Gamma)$ on the one hand, and
the Neumann components $\ctru_{\neu} = (u^{j}_{\neu})_{j=0}^{\mJ}\in\mbH_{\neu}(\Gamma)$ on the other hand.
There is a natural duality between Dirichlet and Neumann multi-trace spaces through the bilinear pairing
\begin{equation}
  \begin{array}{l}
    \llangle \ctru,\ctrp\rrangle:=\sum_{j=0}^{\mJ}\langle u^{j},p^{j}\rangle_{\partial\Omega_{j}}\\[5pt]
    \hspace{1.5cm} \forall \ctru = (u^{0},\dots,u^{\mJ})\in\mbH_{\dir}(\Gamma),\\
    \hspace{1.5cm} \forall \ctrp = (p^{0},\dots,p^{\mJ})\in\mbH_{\neu}(\Gamma).
  \end{array}
\end{equation}
The bilinear pairing defined above \textit{does not} involve any complex
conjugation operation. We shall indifferently write $ \llangle \ctrp, \ctru\rrangle :=
\llangle \ctru,\ctrp\rrangle$ for $\ctru\in\mbH_{\dir}(\Gamma),\ctrp\in\mbH_{\neu}(\Gamma)$.

\quad\\
For the sake of conciseness, we shall denote $\mT_{j}$  instead of
$\mT_{\Omega_{j}}$. The operator $\mT := \mrm{diag}_{j=0\dots \mJ}(\mT_{j}):
\mbH_{\dir}(\Gamma)\to \mbH_{\neu}(\Gamma)$  induces 
a scalar product underlying the norm of $\mbH_{\neu}(\Gamma)$ through
\begin{equation}
  \begin{array}{l}
    (\ctrp,\ctrq)_{\mbH_{\neu}(\Gamma)} = \llangle \mT^{-1}(\ctrp),\overline{\ctrq}\rrangle
    = \sum_{j=0}^{\mJ}\langle\mT_{j}^{-1}(p^{j}),\overline{q}^{j}\rangle_{\partial\Omega_{j}}\\
    \textcolor{white}{(\ctrp,\ctrq)_{\mbH_{\dir}(\Gamma)}} =
    \sum_{j=0}^{\mJ}(p^{j},q^{j})_{\mH^{-1/2}(\partial\Omega_{j})}.
  \end{array}
\end{equation}
for any $\ctrp = (p^{j})_{j=0\dots \mJ}$ and any $\ctrq = (q^{j})_{j=0\dots \mJ}$ in
$\mbH_{\neu}(\Gamma)$. As regards  $\mbH(\Gamma)$, we shall consider a  duality pairing
given by the following skew symetric bilinear form 
\begin{equation}
  \begin{aligned}
    & \llbr \ctru,\ctrv\rrbr:=\lbr\ctru_{0},\ctrv_{0}\rbr_{\partial\Omega_{0}}+\dots+
    \lbr\ctru_{\mJ},\ctrv_{\mJ}\rbr_{\partial\Omega_{\mJ}}\\
    & \text{for}\;\; \ctru = (\ctru_{j})_{j=0}^{\mJ}\;\text{and}\; \ctrv = (\ctrv_{j})_{j=0}^{\mJ}\;\text{in}\;\mbH(\Gamma).
  \end{aligned}
\end{equation}
As regards trace operators, for the sake of conciseness, we shall
denote $\tau^{j}:=\tau^{\Omega_{j}}$ and adopt similar conventions for
$\tau^{j}_{\dir}$ and $\tau^{j}_{\neu}$. We also introduce
global trace operators that map into multi-trace spaces
\begin{equation}
  \begin{array}{l}
    \tau_{\alpha}(u):=(\tau^{0}_{\alpha}(u),\dots,\tau^{\mJ}_{\alpha}(u))
    \quad\text{for}\;\;\alpha = \dir,\neu\\[5pt]
    \tau(u)_{\textcolor{white}{\alpha}}:=(\tau^{0}(u),\dots, \tau^{\mJ}(u)).
  \end{array}
\end{equation}

\section{Transmission  conditions}\label{TransmissionsConditions}
Since we are considering a problem involving transmission
conditions (\ref{IntialP2:2}), it is natural to introduce the subspace of $\mbH(\Gamma)$
consisting in all tuples of traces agreeing with these conditions: this is what
shall be called single-trace spaces defined by
\begin{equation}\label{SingleTraceSpaces}
  \begin{aligned}
    \mbX_{\dir}(\Gamma) & :=\{\;(v_{j})_{j=0}^{\mJ}\in\mbH_{\dir}(\Gamma)\;\vert\;
    \exists\vphi\in\mH^{1}(\RR^{d}),\;v_{j} = \vphi\vert_{\partial\Omega_{j}}\;\forall j\;\}\\
    \mbX_{\neu}(\Gamma) & :=\{\;(q_{j})_{j=0}^{\mJ}\in\mbH_{\neu}(\Gamma)\;\vert\;
    \exists\boldsymbol{\psi}\in\mH(\div,\RR^{d}),\;
    q_{j} = \bn_{j}\cdot\boldsymbol{\psi}\vert_{\partial\Omega_{j}}\;\forall j  \;\}\\
    \mbX(\Gamma)_{\textcolor{white}{\dir}} & :=\{\;\ctru = (\ctru_{\dir},\ctru_{\neu})\in\mbH(\Gamma)\;\vert\;
    \ctru_{\dir}\in\mbX_{\dir}(\Gamma),\; \ctru_{\neu}\in\mbX_{\neu}(\Gamma) \;\}\\[5pt]
  \end{aligned}
\end{equation}
By construction, for a function $u\in \mL^{2}_{\loc}(\RR^{d})$ such that $u\vert_{\Omega_{j}}\in
\mH^{1}_{\loc}(\Delta,\overline{\Omega}_{j})$ for all $j=0\dots \mJ$, the transmission conditions 
(\ref{IntialP2:2}) are equivalent to the statement ``$\tau(u)\in \mbX(\Gamma)$''.
The single-trace space has been extensively studied in the context of multi-trace formulations
\cite{MR3069956}. The following caracterisation of this space was proved in \cite[Prop.6.3]{MR3101780}.

\begin{proposition}\label{polarity}\quad\\
  For any $\ctru\in\mbH(\Gamma)$ we have
  $\ctru\in\mbX(\Gamma)\;\iff\;\llbr\ctru,\ctrv\rrbr = 0\;
  \forall \ctrv\in\mbX(\Gamma)$. 
\end{proposition}
\noindent 
\textbf{Proof:}

From (\ref{SingleTraceSpaces}), it is clear that any $\ctru = (\ctru_{\dir},\ctru_{\neu})
\in\mbH(\Gamma)$ actually belongs to $\mbX(\Gamma)$ if and only if $\ctru_{\dir}
\in\mbX_{\dir}(\Gamma)$ and $\ctru_{\neu}\in\mbX_{\neu}(\Gamma)$. As a consequence,
to prove the lemma, it suffices to show that for any
$\ctru_{\dir}\in\mbH_{\dir}(\Gamma)$ and any $\ctru_{\neu}\in\mbH_{\neu}(\Gamma)$
we have
\begin{itemize}
\item[i)] $\ctru_{\dir}\in\mbX_{\dir}(\Gamma)\iff\llangle\ctru_{\dir},\ctrq\rrangle = 0\;
\forall\ctrq\in\mbX_{\neu}(\Gamma)$\\[-20pt]
\item[ii)] $\ctru_{\neu}\in\mbX_{\neu}(\Gamma)\iff\llangle\ctru_{\neu},\ctrv\rrangle = 0\;
\forall\ctrv\in\mbX_{\dir}(\Gamma)$
\end{itemize}
We will only present the proof of \textit{i)} since the proof for \textit{ii)}
is very similar. Take an arbitrary $\ctru_{\dir} = (u_{\dir}^{j})_{j=0}^{\mJ}
\in\mbH_{\dir}(\Gamma)$. If $\ctru_{\dir}\in\mbX_{\dir}(\Gamma)$, there exists
$\varphi\in\mH^{1}(\RR^{d})$ such that $\varphi\vert_{\partial\Omega_{j}} = u_{\dir}^{j}\,
\forall j=0\dots \mJ$. Then for any $\ctrq = (q^{j})_{j=0}^{\mJ}\in\mbX_{\neu}(\Gamma)$, 
there exists $\bpsi\in \mH(\div,\RR^{d})$ such that $\bn_{j}\cdot\bpsi\vert_{\partial\Omega_{j}} = q^{j}
\,\forall j=0\dots \mJ$. Applying a Green formula in each $\Omega_{j}$ on the one hand, 
and in $\RR^{d}$ on the other hand, we obtain  
\begin{equation}
  \begin{array}{rl}
    \llangle \ctru_{\dir},\ctrq\rrangle 
    & = \sum_{j=0}^{\mJ}\langle u_{\dir}^{j},q^{j}\rangle_{\partial\Omega_{j}} = 
      \sum_{j=0}^{\mJ}\int_{\partial\Omega_{j}}\bn_{j}\cdot\bpsi\varphi \,d\sigma \\
    & = \sum_{j=0}^{\mJ}\int_{\Omega_{j}}\nabla \varphi\cdot\bpsi + \varphi\,\div\bpsi\, d\bx 
      = \int_{\RR^{d}}\nabla \varphi\cdot\bpsi + \varphi\,\div\bpsi\, d\bx =0.
  \end{array}
\end{equation}
Now assume that $\ctru_{\dir} = (u_{\dir}^{0},\dots ,u_{\dir}^{\mJ})
\in\mbH_{\dir}(\Gamma)$ satisfies $\llangle \ctru_{\dir},\ctrq\rrangle = 0\forall\ctrq
\in\mbX_{\neu}(\Gamma)$. For each $j = 0\dots \mJ$, introduce a lifting $v_{j}\in \mH^{1}(\Omega_{j})$ 
such that $v_{j}\vert_{\partial\Omega_{j}} = u_{\dir}^{j}$, and set 
$v(\bx) = 1_{\Omega_{0}}(\bx)v_{0}(\bx)+\dots + 1_{\Omega_{\mJ}}(\bx)v_{\mJ}(\bx)$. We have clearly 
$v\in\mL^{2}(\RR^{d})$ and, to prove that $\ctru_{\dir}\in\mbX_{\dir}(\Gamma)$, it suffices 
to show that $v\in\mH^{1}(\RR^{d})$. Define $\bp\in\mL^{2}(\RR^{d})$ by 
$\bp(\bx) = 1_{\Omega_{0}}(\bx)\nabla v_{0}(\bx)+\dots+1_{\Omega_{\mJ}}(\bx)\nabla v_{\mJ}(\bx)$. 
Pick an arbitrary $\bpsi\in \mH(\div,\RR^{d})$, and set $\ctrq = (q^{j})_{j=0}^{\mJ}$ where 
$q^{j}:=\bn_{j}\cdot\bpsi\vert_{\partial\Omega_{j}}$. Since $\ctrq\in\mbX_{\neu}(\Gamma)$, we have 
\begin{equation}
  \begin{array}{rl}
    \int_{\RR^{d}} v\,\div(\bpsi) d\bx 
    & = \sum_{j=0}^{\mJ}\int_{\Omega_{j}}v\,\div(\bpsi) d\bx\\
    & = \llangle \ctru_{\dir},\ctrq\rrangle-\sum_{j=0}^{\mJ}\int_{\Omega_{j}}\bpsi\cdot\nabla v_{j} d\bx\\
    & = -\int_{\RR^{d}}\bpsi\cdot\bp d\bx
  \end{array}
\end{equation}
Since the above identity holds for any $\bpsi\in \mH(\div,\RR^{d})$, we conclude that $v$
admits a weak gradient over $\RR^{d}$ as a whole with $\bp = \nabla v$ in $\RR^{d}$ and, as
a consequence $v\in\mH^{1}(\RR^{d})$ and $\ctru_{\dir}\in\mbX(\Gamma)$.  \hfill $\Box$

\quad\\
As underlined during its proof, the above caracterisation implies that
$\ctru\in\mbH_{\dir}(\Gamma)$ belongs to $\mbX_{\dir}(\Gamma)$ if and only if
$\llangle \ctru,\ctrp\rrangle = 0\forall \ctrp\in\mbX_{\neu}(\Gamma)$ and that,
similarly, $\ctrp\in\mbH_{\neu}(\Gamma)$ belongs to $\mbX_{\neu}(\Gamma)$ if and
only if $\llangle \ctru,\ctrp\rrangle = 0\, \forall \ctru\in\mbX_{\dir}(\Gamma)$. 

\begin{proposition}\label{OrthogonalDecomposition}\quad\\
  We have the direct sum $\mbH_{\neu}(\Gamma) = \mbX_{\neu}(\Gamma)
  \oplus\mT(\mbX_{\dir}(\Gamma))$ and it is orthogonal with respect
  to the scalar product induced by $\mT^{-1}$.
\end{proposition}
\noindent \textbf{Proof:}

First, according to  Proposition \ref{polarity}, we have 
$(\ctrp,\mT(\ctru))_{\mbH_{\neu}(\Gamma)} = \llangle\ctrp,\ctru\rrangle = 0$ whenever
$\ctrp\in \mbX_{\neu}(\Gamma)$ and $\ctru \in \mT(\mbX_{\dir}(\Gamma))$.
This proves that $\mbX_{\neu}(\Gamma)$ is orthogonal to $\mT(\mbX_{\dir}(\Gamma))$ 
hence $\mbX_{\neu}(\Gamma)\cap \mT(\mbX_{\dir}(\Gamma)) = \{0\}$.

Next pick an arbitrary $\ctrp\in\mbH_{\neu}(\Gamma)$ and, by Riesz representation theorem,
define $\ctru$ as the unique element of $\mbX_{\dir}(\Gamma)$ satisfying
$\llangle \mT(\ctru),\overline{\ctrv}\rrangle  = \llangle\ctrp,\overline{\ctrv}\rrangle$
for all $\ctrv\in\mbX_{\dir}(\Gamma)$. As a consequence $\ctrq = \ctrp - \mT(\ctru)$
satisfies $\llangle \ctrq,\ctrv\rrangle = 0\forall \ctrv\in \mbX_{\dir}(\Gamma)$ and thus 
belongs to $\mbX_{\neu}(\Gamma)$ according to Proposition \ref{polarity}. This shows that
$\mbH_{\neu}(\Gamma) = \mbX_{\neu}(\Gamma) +\mT(\mbX_{\dir}(\Gamma))$. \hfill $\Box$

\section{Potential theory}

The problem (\ref{IntialPb}) primarily considered in the present manuscript does not a priori lend itself
to boundary integral equation techniques simply because (\ref{IntialPb}) is a problem of propagation
in heterogeneous media i.e. the PDEs involve a priori varying coefficients. However several aspects of the
solution strategy we wish to describe involve nonlocal operators. In particular,  we shall need such theoretical
tools for treatment of junctions. As a consequence, we dedicate the present section to recalling a few facts about
boundary integral operators.

\subsection{Layer potentials in a single subdomain}
We first introduce the Green kernel $\mathscr{G}(\bx)$ of the Yukawa's equation i.e.
we define $\green$ as the unique function solving $-\Delta \green + \gamma^{-2}\green = \delta_{0}$ in $\RR^{d}$
and $\lim_{\vert\bx\vert\to \infty}\green(\bx) = 0$, where $\delta_{0}$ is the Dirac measure centered at $\bx = 0$,
and $\gamma>0$ is a parameter that we have fixed once and for all in \S\ref{VolumeSpaces}.
This kernel admits an explicit expression in terms of special functions namely
\begin{equation}
  \begin{aligned}
    & \green(\bx):=\mrm{K}_{0}(\vert\bx\vert/\gamma),& \quad  \bx\in\RR^{2}\setminus\{0\}     & \quad \text{for $d=2$}, \\
    & \green(\bx):=\frac{\exp(-\vert\bx\vert/\gamma)}{4\pi\vert\bx\vert},& \quad  \bx\in\RR^{3}\setminus\{0\} & 
    \quad \text{for $d=3$}.
  \end{aligned}
\end{equation}
where $\mrm{K}_{0}$ refers to the modified Bessel function of the second kind of order $0$ also known
as MacDonald function, see \cite[\S 10.25]{MR2723248}. With this kernel, and for any Lipschitz domain
$\Omega\subset\RR^{d}$ with bounded boundary, we can define single and double layer potentials as
follows: for any $(v,q)\in\mH^{1/2}(\partial\Omega)\times\mH^{-1/2}(\partial\Omega)$ we set 
\begin{equation}
  \begin{array}{l}
    \Psi^{\Omega}(v,q)(\bx):= \Psi^{\Omega}_{\dir}(v)(\bx) + \Psi^{\Omega}_{\neu}(q)(\bx),\\[10pt]
    \text{where}\quad \Psi^{\Omega}_{\dir}(v)(\bx):=\int_{\partial\Omega}\bn_{\Omega}(\by)
    \cdot(\nabla\green)(\bx-\by)v(\by) d\sigma(\by), \\[5pt]
    \textcolor{white}{\text{where}\quad}
    \Psi^{\Omega}_{\neu}(q)(\bx):=\int_{\partial\Omega}\green(\bx-\by)q(\by) d\sigma(\by),\\
  \end{array}
\end{equation}
for all $\bx\in \RR^{d}\setminus\partial\Omega$. For any $\ctrv\in\mH^{1/2}(\partial\Omega)\times\mH^{-1/2}(\partial\Omega)$,
we have $(\gamma^{-2}-\Delta)\Psi^{\Omega}(\ctrv) =0$ both in $\Omega$ and $\RR^{d}\setminus\overline{\Omega}$.
Besides $\Psi^{\Omega}(\ctrv)\vert_{\Omega}\in \mH^{1}(\Delta,\overline{\mathcal{O}})$ for
$\mathcal{O} = \Omega$ or $\mathcal{O} = \RR^{d}\setminus\overline{\Omega}$. For any
$\bx,\by\in\RR^{d},\bx\neq \by$, define $\green_{\bx}:\RR^{d}\setminus\{\bx\}\to \RR_{+}$ by
$\green_{\bx}(\by):=\green(\bx-\by)$. Elementary calculus shows that
$\Psi^{\Omega}(\ctru)(\bx) = \lbr\tau^{\Omega}(\green_{\bx}),\ctru\rbr_{\partial\Omega}$
for all $\ctru\in\mH^{1/2}(\partial\Omega)\times \mH^{-1/2}(\partial\Omega)$ and all
$\bx\in \RR^{d}\setminus\partial\Omega$. The next result, known as representation theorem,
shows that layer potential can be used to reconstruct any solution to the homogeneous Yukawa equation.

\begin{proposition}\label{ReprThm}\quad\\
  For any Lipschitz domain $\Omega\subset \RR^{d}$ with bounded boundary, and any
  function $u\in \mH^{1}(\Omega)$ satisfying $(\gamma^{-2}-\Delta) u = 0$ in $\Omega$, we  have
  $\Psi^{\Omega}(\tau^{\Omega}(u)) = 1_{\Omega}(\bx)u(\bx)\;\forall\bx\in\RR^{d}$.
\end{proposition}

\noindent 
Here $1_{\Omega}(\bx) = 1$ if $\bx\in\Omega$ and $1_{\Omega}(\bx) = 0$ otherwise.
In the representation formula above, the traces of solutions to the homogeneous PDE play
a pivotal role. The potential operators actually provide a Calder\'on projector that
maps onto such a space and can thus be used to caracterise them. 

\begin{proposition}\label{DefCalderonProj}\quad\\
  The operator $\tau^{\Omega}\cdot\Psi^{\Omega}:\mH^{1/2}(\partial\Omega)\times
  \mH^{-1/2}(\partial\Omega)\to \mH^{1/2}(\partial\Omega)\times \mH^{-1/2}(\partial\Omega)$
  is a continuous  projector whose range is the space  
    $\calC_{\mrm{in}}(\Omega):=\{\tau^{\Omega}(u)\;\vert\;u\in\mH^{1}(\Omega),\;
    (\gamma^{-2}-\Delta) u = 0\;\text{in}\;\Omega\;\}$.
\end{proposition}

\subsection{Layer potentials in a multi-domain setting}

Considering $\Omega = \Omega_{j}$ for $j=0\dots \mJ$, the result of the
previous paragraph can be used directly in the multi-domain context.
For the sake of conciseness, in the following, we shall write $\Psi^{j}_{\dir},
\Psi^{j}_{\neu},\Psi^{j}$ instead of $\Psi^{\Omega_{j}}_{\dir},\Psi^{\Omega_{j}}_{\neu},
\Psi^{\Omega_{j}}$.

We now show that an explicitly formula for the orthogonal projector onto $\mbX_{\neu}(\Gamma)$,
can be obtained. We rely on so-called multi-potential operators
$\Psi_{\dir}:\mbH_{\dir}(\Gamma)\to \Pi_{j=0}^{\mJ}\mH^{1}_{\loc}(\Delta,\overline{\Omega}_{j})$ and 
$\Psi_{\neu}:\mbH_{\neu}(\Gamma)\to \Pi_{j=0}^{\mJ}\mH^{1}_{\loc}(\Delta,\overline{\Omega}_{j})$ defined
as folllows: for any $\ctru = (\ctru_{\dir},\ctru_{\neu})\in \mbH(\Gamma)$ we set
\begin{equation}
  \begin{array}{l}
    \Psi(\ctru)(\bx) = \Psi_{\dir}(\ctru_{\dir})(\bx) + \Psi_{\neu}(\ctru_{\neu})(\bx)\\[10pt]
    \text{where}\quad \Psi_{\dir}(\ctru_{\dir})(\bx):=\sum_{j=0}^{\mJ}
    \Psi^{j}_{\dir}(u^{j}_{\dir})(\bx),\\
    \textcolor{white}{\text{where}\quad}
    \Psi_{\neu}(\ctru_{\neu})(\bx):=\sum_{j=0}^{\mJ}\Psi^{j}_{\neu}(u^{j}_{\neu})(\bx).
  \end{array}
\end{equation}
for any $\bx\in\RR^{d}\setminus\Gamma$. Such operators have been first
considered in the context of the integral formulation of the second kind
introduced in \cite{C11_144}, see also \cite{MR3313601,MR3725822,MR3720391,MR3403719}.
The multi-potential operators satisfy many non-trivial properties. To begin with,
the next proposition shows that they are closely related to global
Dirichlet-to-Neumann maps.

\begin{lemma}\label{RadiationReprFormula}\quad\\
  We have $\tau\cdot\Psi(\ctru) = \ctru$  for all
  $\ctru = (\ctru_{\dir},\ctru_{\neu})\in \mbH(\Gamma)$ satisfying 
  $\ctru_{\neu} = \mT(\ctru_{\dir})$.
\end{lemma}
\noindent \textbf{Proof:}

Pick a  $\ctru = (\ctru_{\dir},\ctru_{\neu}) =
(\ctru_{\dir}^{j},\ctru_{\neu}^{j})_{j=0}^{\mJ} \in \mbH(\Gamma)$ with
$\ctru_{\neu} = \mT(\ctru_{\dir})$. We have $\ctru^{j} := (\ctru_{\dir}^{j},\ctru_{\neu}^{j}) =
(\ctru_{\dir}^{j},\mT_{j}(\ctru_{\dir}^{j}))\in \mathcal{C}_{\mrm{in}}(\Omega_{j})$
for each $j = 0\dots \mJ$. As a consequence, applying Proposition \ref{ReprThm}, we obtain
$\tau^{k}\Psi^{j}(\ctru^{j}) = \delta_{j,k}\ctru^{j}$
for any $j,k = 0\dots \mJ$. Summing the latter identity over $j$
yields $\tau^{k}\Psi(\ctru) = \ctru^{k}$ for all $k = 0\dots \mJ$,
which concludes the proof.\hfill $\Box$

\begin{lemma}\label{KernelMultiPot}\quad\\
  We have $\Psi(\ctru) = 0\,\forall\ctru\in\mbX(\Gamma)$.
\end{lemma}
\noindent \textbf{Proof:}

Denoting as before $\green_{\bx}(\by):=\green(\bx-\by)$, recall that we have
$\Psi^{\Omega}(\ctru)(\bx) = \lbr\tau^{\Omega}(\green_{\bx}),\ctru\rbr_{\partial\Omega}$
$\ctru\in\mH^{1/2}(\partial\Omega)\times \mH^{-1/2}(\partial\Omega)$ and all
$\bx\in \RR^{d}\setminus\partial\Omega$. Plugging this expression
into the definition of the multi-potential operator yields 
 $ \Psi(\ctru)(\bx) = \llbr\tau(\green_{\bx}),\ctru\rrbr
  \quad\forall\ctru\in \mbH(\Gamma), \;\;\forall\bx\in \RR^{d}\setminus\Gamma$.
Now observe that for any $\bx\in \RR^{d}\setminus\Gamma$
we have $\tau(\green_{\bx})\in \mbX(\Gamma)$ hence  
applying Proposition \ref{polarity} concludes the proof.
\hfill $\Box$

\quad\\
A direct consequence of the lemma above is that $\Psi_{\dir}(\ctru_{\dir}) = 0$
for all $\ctru_{\dir}\in \mbX_{\dir}(\Gamma)$, and $\Psi_{\neu}(\ctru_{\neu}) = 0$
for all $\ctru_{\neu}\in \mbX_{\neu}(\Gamma)$. We deduce in particular
that $\mbX_{\neu}(\Gamma)\subset \mrm{Ker}(\tau_{\neu}\cdot\Psi_{\neu})$. 

\begin{lemma}\quad\\
  We have $\ctrp - \tau_{\neu}\cdot\Psi_{\neu}(\ctrp)\in
  \mbX_{\neu}(\Gamma)$ for any $\ctrp\in \mbH_{\neu}(\Gamma)$.
\end{lemma}
\noindent \textbf{Proof:}

Pick an arbitrary $\ctrp\in \mbH_{\neu}(\Gamma)$ and, applying Proposition \ref{OrthogonalDecomposition},
decompose it as $\ctrp = \ctrv_{\neu} + \mT(\ctru_{\dir})$ where $\ctru_{\dir}
\in\mbX_{\dir}(\Gamma)$ and $\ctrv_{\neu}\in\mbX_{\neu}(\Gamma)$.
According to Lemma \ref{KernelMultiPot} we have $\Psi_{\dir}(\ctru_{\dir}) = 0$
so that, setting $\ctru:=(\ctru_{\dir},\mT(\ctru_{\dir}))$, we have
$\ctrp - \tau_{\neu}\cdot\Psi_{\neu}(\ctrp) =\ctrv_{\neu}-\tau_{\neu}\cdot\Psi_{\neu}(\ctrv_{\neu}) +\mT(\ctru_{\dir})-\tau_{\neu}\cdot\Psi(\ctru)$.
Applying Lemma \ref{RadiationReprFormula}
yields $\mT(\ctru_{\dir})-\tau_{\neu}\cdot\Psi(\ctru) = 0$. Besides we have
$\Psi_{\neu}(\ctrv_{\neu}) = 0$ according to Lemma \ref{KernelMultiPot} since
$\ctrv_{\neu}\in\mbX_{\neu}(\Gamma)$. To summarise, we have just established
$\ctrp - \tau_{\neu}\cdot\Psi_{\neu}(\ctrp) =\ctrv_{\neu}\in\mbX_{\neu}(\Gamma)$,
which concludes the proof. \hfill $\Box$

\quad\\
Combining the previous two lemmas, we see that
$(\tau_{\neu}\cdot\Psi_{\neu})(\Id - \tau_{\neu}\cdot\Psi_{\neu}) = 0$.
From this we deduce immediately the following proposition.

\begin{proposition}\label{ProjectorKern}\quad\\
  We have $\mrm{Ker}(\tau_{\neu}\cdot\Psi_{\neu}) = \mrm{Range}(\mrm{Id}-\tau_{\neu}\cdot\Psi_{\neu})
  = \mbX_{\neu}(\Gamma)$, and $\tau_{\neu}\cdot\Psi_{\neu}:
  \mbH_{\neu}(\Gamma)\to \mbH_{\neu}(\Gamma)$ is a continuous projector.
\end{proposition}

\noindent 
The next result gives further details about the image of this projector. 

\begin{lemma}\label{OrthogonalityProj}\quad\\
  We have $\mrm{Range}(\tau_{\neu}\cdot\Psi_{\neu}) = \mT(\mbX_{\dir}(\Gamma))$ so that
  $\tau_{\neu}\cdot\Psi_{\neu}$ is an orthogonal projector with respect to the scalar
  product induced by $\mT^{-1}$ over $\mbH_{\neu}(\Gamma)$.
\end{lemma}
\noindent \textbf{Proof:}

Taking account of both Proposition \ref{OrthogonalDecomposition} and \ref{ProjectorKern},
we see that it suffices to prove $\tau_{\neu}\cdot\Psi_{\neu}(\mT(\ctru)) = \mT(\ctru)$
for all $\ctru\in \mbX_{\dir}(\Gamma)$. Hence consider any $\ctru = (u_{j})_{j=0}^{\mJ}\in\mbX_{\dir}(\Gamma)$.
According to Lemma \ref{KernelMultiPot} we have $\tau_{\neu}\cdot\Psi_{\dir}(\ctru) = 0$.
As a consequence, applying Corollary \ref{RadiationReprFormula},  we obtain 
\begin{equation}
  \begin{array}{rl}
    \tau^{k}_{\neu}\cdot\Psi_{\neu}(\mT(\ctru))
    &  = \tau^{k}_{\neu}\cdot (\;\Psi_{\dir}(\ctru) + \Psi_{\neu}(\mT(\ctru))\;)\\
    &  = \tau^{k}_{\neu}\cdot\sum_{j=0}^{\mJ}\Psi_{\dir}^{j}(u_{j}) + 
      \Psi_{\neu}^{j}(\mT_{j}(u_{j}))\\
    & = \sum_{j=0}^{\mJ}\tau^{k}_{\neu}\cdot\Psi_{\dir}^{j}(u_{j}) +
      \tau^{k}_{\neu}\cdot\Psi_{\neu}^{j}(\mT_{j}(u_{j})) = \mT_{k}(u_{k})
  \end{array}
\end{equation}
for any $k=0\dots \mJ$. Since this holds for all $k$, we obtain that
$\tau_{\neu}\cdot\Psi_{\neu}(\mT(\ctru)) = \mT(\ctru)$,
which concludes the proof. \hfill $\Box$

\quad\\
From the previous results, we immediately obtain an estimate on the norm of the projection,
which  will be key in the analysis of Section \ref{ReformulationGlobalScat}. 

\begin{corollary}\label{ContractionProj}\quad\\
  Define $\Pi := \mrm{Id}-2\tau_{\neu}\cdot\Psi_{\neu}$.
  Then we have $\Pi^{2} = \mrm{Id}$ and the operators $(\Id\pm\Pi)/2$
  are continuous projectors with $\mbX_{\neu}(\Gamma):=\mrm{Ker}(\mrm{Id}-\Pi)$
  and $\mT(\mbX_{\dir}(\Gamma)):=\mrm{Ker}(\mrm{Id}+\Pi)$. 
  Besides the following continuity estimate holds:
  \begin{equation*}
    \Vert \Pi(\ctrp)\Vert_{\mbH_{\neu}(\Gamma)}=
    \Vert \ctrp\Vert_{\mbH_{\neu}(\Gamma)}\quad
    \forall \ctrp\in\mbH_{\neu}(\Gamma).
  \end{equation*}
\end{corollary}

\noindent 
In the subsequent analysis, this projector will be the key tool
for caracterising elements of $\mbX(\Gamma)$ and thus enforcing
transmission conditions across interfaces. The next result indeed
provides a caracterisation of the single trace space.

\begin{proposition}\label{CaracTransmissionRobinTraces}\quad\\
  Consider any $\omega >0$. With the notations of the previous corollary,
  for any $\ctru = (\ctru_{\dir},\ctru_{\neu})\in \mbH(\Gamma)$, we
  have $\ctru\in \mbX(\Gamma)$ if and only if
  $\ctru_{\neu}-\imath\omega \mT(\ctru_{\dir}) =
  \Pi(\ctru_{\neu}+\imath\omega\mT(\ctru_{\dir}))$.
\end{proposition}
\noindent \textbf{Proof:}

According to Corollary \ref{ContractionProj}, for $\ctru =
(\ctru_{\dir},\ctru_{\neu})\in\mbH(\Gamma)$, we have $\ctru_{\neu}\in\mbX_{\neu}(\Gamma)
\iff (\Id-\Pi)\ctru_{\neu} = 0$ and $\ctru_{\dir}\in\mbX_{\dir}(\Gamma) \iff
(\Id+\Pi)\mT(\ctru_{\dir}) = 0$. On the other hand,
$\mrm{Range}(\Id+\Pi)\cap \mrm{Range}(\Id-\Pi) = \{0\}$
since $(\Id+\Pi)/2$ is a projector, which leads to
$\ctru\in\mbX(\Gamma)\iff (\Id-\Pi)\ctru_{\neu} =
\imath \omega (\Id+\Pi)\mT(\ctru_{\dir})$. Rearranging this latter identity 
yields the conclusion of the proof. \hfill $\Box$

\section{Reformulation of wave equations}

In this section we focus on the wave equations (\ref{IntialP2:1}) that we will
reformulate in terms of traces only. We adopt the approach developped by
Collino, Ghanemi and Joly in \cite{MR1764190} and further studied and extended in
\cite{LECOUVEZ2014403,LecouvezThesis}. This approach generalises the original work of Després
\cite{MR1227838,MR1291197,MR1105979,MR1071633} on Optimised Schwarz Method for Helmholtz equation.
In the present section, we will derive a convenient caracterisation of
\begin{equation}
  \begin{aligned}
    & \mathscr{C}^{+}(\Gamma):= \mathscr{C}^{+}(\Omega_{0})\times \dots \times \mathscr{C}^{+}(\Omega_{\mJ})
    \quad\text{where}\\
    & \mathscr{C}^{+}(\Omega_{j}):=\{\;(\tau_{\dir}^{j}(\varphi),\mu_{j}\tau_{\neu}^{j}(\varphi)) \in 
    \mH^{1/2}(\partial\Omega_{j})\times \mH^{-1/2}(\partial\Omega_{j}),\\
    &\textcolor{white}{\mathscr{C}^{+}(\Omega_{j}):=\{\;\hspace{2cm}}
    \div(\mu\nabla\varphi) + \kappa^{2}\varphi = 0\;\textrm{in}\;\Omega_{j}\;\text{and}\\
    & \textcolor{white}{\mathscr{C}^{+}(\Omega_{j}):=\{\;\hspace{2cm}}
    \hspace{1cm}\varphi\;\kappa_{0}-\textrm{outgoing if $j=0$.}\;\}.
  \end{aligned}
\end{equation}
The space $\mathscr{C}^{+}(\Omega_{j})$ is closed in $\mH^{1/2}(\partial\Omega_{j})\times \mH^{-1/2}(\partial\Omega_{j})$
and we will use these spaces to reformulate the wave equation in each subdomain. We have the following important 
decomposition of the multi-trace space. 

\begin{proposition}\label{DecompSTFPropagCauchy}\quad\\
  We have the direct sum $\mbH(\Gamma) = \mbX(\Gamma)\oplus\mathscr{C}^{+}(\Gamma)$.
\end{proposition}
\noindent \textbf{Proof:}

Let us first show that $\mbX(\Gamma)\cap \mathscr{C}^{+}(\Gamma) = \{0\}$.
Pick some $\ctru\in \mbX(\Gamma)\cap \mathscr{C}^{+}(\Gamma)$
decomposed in Dirichlet/Neumann components  $\ctru = (\ctru_{\dir},\ctru_{\neu})$ with
$\ctru_{\dir} = (u_{\dir}^{j})_{j=0}^{\mJ}\in\mbH_{\dir}(\Gamma)$ and 
$\ctru_{\neu} = (u_{\neu}^{j})_{j=0}^{\mJ}\in\mbH_{\neu}(\Gamma)$.
For each $j=0\dots \mJ$, let $\phi_{j}\in \mH^{1}_{\loc}(\Omega_{j})$ refer to the unique functions satisfying
\begin{equation}\label{PbDecompDirectSum3}
  \begin{aligned}
    & \div(\mu\nabla\phi_{j})+\kappa^{2}\phi_{j}=0\quad \text{in}\;\Omega_{j},\\
    & \phi_{0}\;\text{is $\kappa_{0}$-outgoing},\\
    & (\tau_{\dir}^{j}(\phi_{j}),\mu_{j}\tau_{\neu}^{j}(\phi)) = (u_{\dir}^{j},u_{\neu}^{j})\;\;\text{on}\;\partial\Omega_{j}.
  \end{aligned}
\end{equation}
Set $\phi := 1_{\Omega_{0}}\phi_{0}+ \dots + 1_{\Omega_{\mJ}}\phi_{\mJ}$, so that
$\div(\mu\nabla\phi) + \kappa^{2}\phi = 0$ in each $\Omega_{j}$ and, since
$\ctru = (\tau_{\dir}^{j}(\phi),\mu_{j}\tau_{\neu}^{j}(\phi) )_{j=0\dots \mJ}\in \mbX(\Gamma)$
the function $\phi$ satisfies transmission conditions across $\Gamma$, so that
$\div(\mu\nabla \phi) + \kappa^{2}\phi = 0$ in $\RR^{d}$ and $\phi$
is $\kappa_{0}$-outgoing. Well-posedness of the Helmholtz equation with outgoing radiation condition
leads to $\phi = 0$, hence $\ctru = 0$, which proves that
\begin{equation}\label{NullIntersection}
\mbX(\Gamma)\cap \mathscr{C}^{+}(\Gamma) = \{0\}.
\end{equation}
Now let us consider the general case of an arbitrary $\ctru\in \mbH(\Gamma)$.
Consider any lifting function $\psi'\in \mL^{2}(\RR^{d})$ with compact support such that
$\psi'\vert_{\Omega_{j}}\in \mH^{1}(\Omega_{j})$ and $\tau^{j}_{\dir}(\psi') = u_{\dir}^{j}$ for all $j=0\dots \mJ$. 
Next define $\psi\in \mH^{1}_{\loc}(\RR^{d})$ as the unique element of $\mH^{1}_{\loc}(\RR^{d})$ satisfying
\begin{equation}\label{PbDecompDirectSum}
  \begin{array}{l}
    \sum_{j=0}^{\mJ}\int_{\Omega_{j}}\mu\nabla (\psi+\psi')\cdot\nabla \vphi
    -\kappa^{2}(\psi+\psi')\vphi\;d\bx\\[5pt]
    \hspace{3.5cm}=\llangle \ctru_{\neu},\tau_{\dir}(\varphi)\rrangle
  \quad\forall \varphi\in \mH^{1}_{\mrm{comp}}(\RR^{d}) \\[7.5pt]
  \text{and}\quad \lim_{\rho\to \infty}\int_{\partial\mB_{\rho}}\vert \partial_{\rho}\psi -
  \imath\kappa_{0}\psi\vert^{2}d\sigma_{\rho}=0
  \end{array}
\end{equation}
where $\mH^{1}_{\mrm{comp}}(\RR^{d})$ refers to the elements of $\mH^{1}(\RR^{d})$
that are boundedly supported. Existence and uniqueness of such a $\psi$ stems
from well posedness of Helmholtz problems in unbounded heterogeneous media,
see e.g. \cite[Chap.3]{MR2986407}.  Applying a Green formula in each $\Omega_{j}$,
we obtain
\begin{equation*}
  \begin{aligned}
    & \div(\mu \nabla(\psi+\psi')) + \kappa^{2}(\psi+\psi') = 0\quad \text{in each}\;\Omega_{j}, j=0\dots \mJ\\
    & \lim_{\rho\to \infty}\int_{\partial\mB_{\rho}}\vert \partial_{\rho}\psi -
    \imath\kappa_{0}\psi\vert^{2}d\sigma_{\rho}=0.    
  \end{aligned}
\end{equation*}
Setting $\ctrv = (\tau_{\dir}^{j}(\psi+\psi'),\mu_{j}\tau_{\neu}^{j}(\psi+\psi'))_{j=0,\dots,\mJ}$, the 
equations above imply that $\ctrv\in\mathscr{C}^{+}(\Gamma)$. Decomposing in Dirichlet/Neumann
contributions $\ctrv = (\ctrv_{\dir},\ctrv_{\neu})$, we have $\ctrv_{\dir} -\ctru_{\dir}=
(\tau^{j}_{\dir}(\psi))_{j=0}^{\mJ}\in \mbX_{\dir}(\Gamma)$ since $\psi\in\mH^{1}_{\loc}(\RR^{d})$. 
Moreover, applying Green formulas once more in (\ref{PbDecompDirectSum}), we see 
$\llangle \ctrv_{\neu},\tau_{\dir}(\varphi)\rrangle = \llangle
\ctru_{\neu},\tau_{\dir}(\varphi)\rrangle$ for all $\varphi\in \mH^{1}(\RR^{d})$.
Using the weak caracterisation of single trace spaces given by Proposition \ref{polarity},
we conclude that $\ctru_{\dir}-\ctrv_{\dir}\in \mbX_{\dir}(\Gamma)$ and
$\ctru_{\neu}-\ctrv_{\neu}\in \mbX_{\neu}(\Gamma)$ hence, setting $\ctrw:=\ctru-\ctrv \in \mbX(\Gamma)$,
so that, with the decomposition $\ctru = \ctrv+\ctrw$, we have established $\mbH(\Gamma) =
\mbX(\Gamma)+\mathscr{C}^{+}(\Gamma)$ which, together with (\ref{NullIntersection}), 
concludes the proof.  \hfill $\Box$

\quad\\
The previous result can be regarded as analogous to Proposition \ref{OrthogonalDecomposition} 
although, in the result above, the direct sum is a priori not orthogonal.  The next property 
relates to energy conservation considerations and will thus play a key role in the 
forthcoming convergence analysis.

\begin{lemma}\label{EnergyFlux}\quad\\
  We have $\imath\lbr\ctru,\overline{\ctru}\rbr_{\partial\Omega_{j}} \leq 0,\;\forall\ctru\in
  \mathscr{C}^{+}(\Omega_{j})\;\forall j=0\dots\mJ$, and thus $\imath \llbr \ctru,
  \overline{\ctru}\rrbr\leq 0\; \forall \ctru\in\mathscr{C}^{+}(\Gamma)$.
\end{lemma}
\noindent \textbf{Proof:}

For any $\ctru\in\mH^{1/2}(\partial\Omega_{j})\times \mH^{-1/2}(\partial\Omega_{j})$, let $\varphi\in\mH^{1}_{\loc}(\overline{\Omega}_{j})$
satisfy $\div(\mu\nabla \vphi)+\kappa^{2}\vphi = 0$ in $\Omega_{j}$  and $(\tau_{\dir}^{j}(\vphi),
\mu_{j}\tau_{\neu}^{j}(\vphi)) = \ctru$ on $\partial\Omega_{j}$. For all $j = 0\dots \mJ$, we have
$\imath\lbr\ctru,\overline{\ctru}\rbr_{\partial\Omega_{j}} =
2\Im m\{\int_{\partial\Omega_{j}}\mu_{j}\tau^{j}_{\neu}(\vphi)\tau^{j}_{\dir}(\overline{\vphi}) d\sigma\}$.
In the case where $j\neq 0$, the domain $\Omega_{j}$ is bounded so that we can apply a simple
Green formula on the later identity,
\begin{equation*}
  \begin{aligned}
    \imath\lbr\ctru,\overline{\ctru}\rbr_{\partial\Omega_{j}}
    & = 2\Im m\{\int_{\partial\Omega_{j}}\mu_{j}\tau^{j}_{\neu}(\vphi)\tau^{j}_{\dir}(\overline{\vphi}) d\sigma\}
    = 2\Im m\{\int_{\Omega_{j}}\overline{\vphi}\div(\mu\nabla\vphi) + \mu\vert \nabla\vphi\vert^{2} d\bx\}\\
    & = 2\Im m\{\int_{\Omega_{j}}\mu\vert \nabla\vphi\vert^{2} - \kappa^{2}\vert\vphi\vert^{2} d\bx\}
    = - 2\int_{\Omega_{j}} \Im m\{\kappa^{2}\}\vert\vphi\vert^{2} d\bx\leq 0.
  \end{aligned}
\end{equation*}
In the case of $\Omega_{0}$ take any radius $\rho_{0}>0$ large enough
to guarantee $\RR^{d}\setminus\Omega_{0}\subset \mB_{\rho_{0}}$. We can apply the
same calculus as above, considering $\mB_{\rho}\cap \Omega_{0}$ instead of
$\Omega_{0}$. Taking account of the radiation condition satisfied 
by $\vphi(\bx)$ for $\vert\bx\vert\to \infty$, and the fact that
$\Im m\{\kappa^{2}\}$ is boundedly supported (since $\kappa(\bx) = \kappa_{0}$
for $\vert\bx\vert>\rho_{0}$), we obtain
\begin{equation*}
  \begin{aligned}
    \imath\lbr\ctru,\overline{\ctru}\rbr_{\partial\Omega_{j}}
    & = -2\int_{\Omega_{0}} \Im m\{\kappa^{2}\}\vert\vphi\vert^{2} d\bx
    +2\Im m\{\int_{\partial\mB_{\rho}}\vphi\partial_{\rho}\overline{\vphi} \,d\sigma\}\\
    & \leq 2\Im m\{\int_{\partial\mB_{\rho}}\vphi\partial_{\rho}\overline{\vphi} \,d\sigma\} =
    -\frac{1}{\kappa_{0}}\int_{\partial\mB_{\rho}}2\Re e\{\imath\kappa_{0}\vphi
    \partial_{\rho}\overline{\vphi} \}\,d\sigma\\
    & = \frac{1}{\kappa_{0}}\int_{\partial\mB_{\rho}}\vert
    \partial_{\rho}\vphi - \imath\kappa_{0}\vphi\vert^{2} d\sigma
    - \frac{1}{\kappa_{0}}\int_{\partial\mB_{\rho}}\vert \partial_{\rho}\vphi\vert^{2}d\sigma
    - \kappa_{0}\int_{\partial\mB_{\rho}}\vert \vphi\vert^{2}d\sigma\\
    & \leq \frac{1}{\kappa_{0}}\int_{\partial\mB_{\rho}}\vert
    \partial_{\rho}\vphi - \imath\kappa_{0}\vphi\vert^{2} d\sigma\quad \forall \rho>\rho_{0}\\
    & \leq \liminf_{\rho\to \infty}\frac{1}{\kappa_{0}}\int_{\partial\mB_{\rho}}\vert
    \partial_{\rho}\vphi - \imath\kappa_{0}\vphi\vert^{2} d\sigma = 0.
  \end{aligned}
\end{equation*}
\hfill $\Box$

\subsection{Robin trace operators}
The caracterisation of $\mbX(\Gamma)$ provided by Proposition
\ref{CaracTransmissionRobinTraces} involved specific combinations 
of Neumann and Dirichlet trace operators. 
Let us bring the attention of the reader to the following
elementary identity: for any $\ctrv = (\ctrv_{\dir},\ctrv_{\neu})\in\mbH(\Gamma)$,
and any $\omega >0$ we have 
\begin{equation}\label{EstimateRobinTrace}
  \begin{aligned}
    \Vert \ctrv_{\neu}+\imath\alpha\mT(\ctrv_{\dir})\Vert_{\mbH_{\neu}(\Gamma)}^{2}
    & = \Vert \ctrv_{\neu}\Vert_{\mbH_{\neu}(\Gamma)}^{2} +
    \omega^{2}\Vert \ctrv_{\dir}\Vert_{\mbH_{\dir}(\Gamma)}^{2} +
    2\alpha\Re e\{\imath\llangle \ctrv_{\dir},\overline{\ctrv}_{\neu}\rrangle\}\\
    & = \Vert \ctrv_{\neu}\Vert_{\mbH_{\neu}(\Gamma)}^{2} +
    \omega^{2}\Vert \ctrv_{\dir}\Vert_{\mbH_{\dir}(\Gamma)}^{2}-2\alpha\Im m\{\llangle \ctrv_{\dir},
    \overline{\ctrv}_{\neu}\rrangle\}\\
    & = \Vert \ctrv_{\neu}\Vert_{\mbH_{\neu}(\Gamma)}^{2} +
    \omega^{2}\Vert \ctrv_{\dir}\Vert_{\mbH_{\dir}(\Gamma)}^{2}
    +\imath\alpha\llbr \ctrv,\overline{\ctrv}\rrbr
    \quad\quad\text{for}\;\;\alpha=\pm \omega.    
\end{aligned}
\end{equation}
We shall assume that the scalar coefficient $\omega>0$, usually referred to as impedance,
is fixed until the end of this article. From the above identity we deduce an expression for
the difference between ingoing and outgoing traces.

\begin{corollary}\label{CorEstimateRobinTrace}\quad\\
  We have $\Vert \ctrv_{\neu}+\imath\omega \mT(\ctrv_{\dir})\Vert_{\mbH_{\neu}(\Gamma)}^{2} -
  \Vert \ctrv_{\neu}-\imath\omega\mT(\ctrv_{\dir})\Vert_{\mbH_{\neu}(\Gamma)}^{2} =
  2\imath \omega\llbr \ctrv,\overline{\ctrv}\rrbr$ for all $\ctrv = (\ctrv_{\dir},\ctrv_{\neu})\in\mbH(\Gamma)$.
\end{corollary}

\noindent
So-called ingoing/outgoing Robin trace operators also play an important role in scattering theory
so, in the present paragraph, we study these trace operators in more detail.
Define $\tau^{j}_{\pm}:\mH^{1}(\Delta,\overline{\Omega}_{j})\to \mH^{-1/2}(\partial\Omega_{j})$ by
\begin{equation}
  \begin{aligned}
    & \tau^{j}_{\pm}(\phi):=\mu_{j}\tau_{\neu}^{j}(\phi)\pm\imath\omega \mT_{j}(\tau_{\dir}^{j}(\phi))
    \quad\text{for}\;\phi\in\mH^{1}(\Delta,\overline{\Omega}_{j}),\\
    & \tau_{\pm}:=\mrm{diag}_{j=0\dots \mJ}(\tau^{j}_{\pm}).\\
  \end{aligned}
\end{equation}
The Robin trace operators can be considered for prescribing boundary data for the solution of wave
equations in each subdomain. Due to the positivity of the DtN maps $\mT_{j}$, the associated
boundary value problems are systematically well posed.

\begin{lemma}\label{WellPosedScatteringPb}\quad\\
  For any $g\in\mL^{2}(\Omega_{j})$ with bounded support, and any $h\in\mH^{-1/2}(\partial\Omega_{j})$, there exists
  a unique $\phi\in\mH^{1}_{\loc}(\overline{\Omega}_{j})$ such that $\div(\mu\nabla \phi) + \kappa^{2}\phi = g$ in $\Omega_{j}$,
  and  $\tau^{j}_{-}(\phi)=h$ on $\partial\Omega_{j}$ (and $\phi$ is $\kappa_{0}$-outgoing if $j=0$). 
\end{lemma}

\noindent 
The proof of the previous lemma is a basic exercise on variationnal formulations, so it is left
to the reader. We need to introduce resolvent operators that solve Helmholtz equation in each
subdomain with a prescribed outgoing Robin boundary trace, the operator
$\mS^{j}:\mH^{-1/2}(\partial\Omega_{j})\to \mH^{-1/2}(\partial\Omega_{j})$ defined by
\begin{equation}
  \begin{aligned}
    \mS^{j}( \tau^{j}_{-}(\phi) ) = \tau^{j}_{+}(\phi)\quad
    &  \text{for all}\; \phi\in \mH^{1}_{\loc}(\overline{\Omega}_{j})\;
    \text{satisfying}\\
    &  \div(\mu\nabla \phi)+\kappa^{2}\phi = 0\quad\text{in}\;\Omega_{j},\\
    &  \phi\;\kappa_{0}-\text{outgoing radiating for $j = 0$}.\\[10pt]
  \end{aligned}
\end{equation}

\begin{proposition}\label{Contractivity}\quad\\
  The operator $\mS = \mrm{diag}_{j=0\dots \mJ}(\mS^{j})$
  continuously maps $\mbH_{\neu}(\Gamma)$ into  $\mbH_{\neu}(\Gamma)$
  and is contractive: for all $\ctrp\in \mbH_{\neu}(\Gamma)$ we have
  \begin{equation*}
    \Vert \mS(\ctrp)\Vert_{\mbH_{\neu}(\Gamma)}\leq \Vert \ctrp\Vert_{\mbH_{\neu}(\Gamma)}.
  \end{equation*}
\end{proposition}
\noindent \textbf{Proof:}

Pick an arbitrary $\ctrp = (p^{j})_{j=0}^{\mJ}\in \mbH_{\neu}(\Gamma)$.
Applying Lemma \ref{WellPosedScatteringPb}, there exist functions
$\phi_{j}\in \mH^{1}_{\loc}(\overline{\Omega}_{j})$ such that $\div(\mu\nabla\phi_{j}) + \kappa^{2}\phi_{j} = 0$
in $\Omega_{j}$, and  $\tau^{j}_{-}(\phi_{j})=p^{j}$ on $\partial\Omega_{j}$ (and $\phi_{j}$ is
$\kappa_{0}$-outgoing if $j=0$). Set $\ctrv = (\ctrv_{\dir},\ctrv_{\neu}):=(\tau^{j}_{\dir}(\phi_{j}),
\mu_{j}\tau^{j}_{\neu}(\phi_{j}))_{j=0,\dots,\mJ}$, we have $\ctrv_{\neu}-\imath\omega \mT(\ctrv_{\dir}) = \ctrp$
and $\ctrv_{\neu}+\imath\omega \mT(\ctrv_{\dir}) = \mS(\ctrp)$.
Since $\ctrv\in \mathscr{C}^{+}(\Gamma)$ by construction, combining
Corollary \ref{CorEstimateRobinTrace} and Lemma \ref{EnergyFlux} concludes
the proof. \hfill $\Box$

\quad\\
The previous result shows that the scattering operator $\mS$ is a contraction
but it is not a priori an isometry. In the context of Problem (\ref{IntialPb}),
this is due to energy loss through radiation of waves toward infinity and absorption
properties of the propagation medium (positive imaginary part of $\kappa^{2}$).

\section{Reformulation of the scattering problem}\label{ReformulationGlobalScat}

In the present section we describe a reformulation of the scattering problem (\ref{IntialPb})
as an equivalently well posed problem.

\subsection{Derivation of the formulation}
To take account of the right hand side $f\in\mH^{1}_{\loc}(\RR^{d})'$,
we introduce the offset function $\phi_{f}\in \mL^{2}_{\loc}(\RR^{d})$
whose restriction to each subdomain  $\phi_{f}\vert_{\Omega_{j}}$ belongs to
$\mH^{1}_{\loc}(\overline{\Omega}_{j})$  and is the unique solution to 
\begin{equation}
  \begin{aligned}
    & \div(\mu\nabla\phi_{f}) + \kappa^{2}\phi_{f} = -f\quad\text{in}\;\Omega_{j},\\
    & \text{$\phi_{f}$ is $\kappa_{0}$-outgoing},\\
    & \tau^{j}_{-}(\phi_{f}) = 0.
  \end{aligned}
\end{equation}
Next, if $u\in\mH^{1}_{\loc}(\RR^{d})$ refers to the unique solution to (\ref{IntialPb})
then $(\tau^{j}_{\dir}(u),\mu_{j}\tau^{j}_{\neu}(u))_{j=0\dots\mJ} \in\mbX(\Gamma)$ so, according to
Proposition \ref{CaracTransmissionRobinTraces}, we have $\tau_{-}(u) = \Pi(\tau_{+}(u))$.
In addition, the function $u-\phi_{f}$ solves an homogenous Helmholtz equation in each
subdomain i.e. $(\div(\mu\nabla\,\cdot\,)+\kappa^{2})(u-\phi_{f}) = 0$ in $\Omega_{j}$ for each $j=\dots \mJ$
and $u-\phi_{f}$ is $\kappa_{0}$-outgoing radiating, so $(\tau^{j}_{\dir}(u-\phi_{f}),
\mu_{j}\tau^{j}_{\neu}(u-\phi_{f})  )_{j=0\dots \mJ}\in\mathscr{C}^{+}(\Gamma)$. As a consequence we have
$\tau_{+}(u-\phi_{f}) = \mS\cdot\tau_{-}(u-\phi_{f}) =\mS\cdot\tau_{-}(u)$.
Thus we conclude that  $\tau_{-}(u) = \Pi\mS(\tau_{-}(u)) + \Pi\tau_{+}(\phi_{f})$.
From this discussion we obtain a reformulation of  our initial scattering
problem  (\ref{IntialPb}),
\begin{equation}\label{OSMFormulation}
  \begin{aligned}
    & \ctrp=\tau_{-}(u)\in \mbH_{\neu}(\Gamma)\quad\text{and}\\
    & \ctrp - (\Pi\cdot\mS)\ctrp = \ctrf\\[5pt]
    & \text{where}\quad \ctrf:= \Pi(\tau_{+}(\phi_{f})).
  \end{aligned}
\end{equation}
The structure of this new formulation is strikingly close to 
standard Optimised Schwarz Methods (OSM). This appears clearly
when comparing (\ref{OSMFormulation}) with \S 2 in \cite{MR1764190},
see in particular Formula (45) and (51) of this reference.

Here also (\ref{OSMFormulation}) appears adapted to domain
decomposition. In the operator $\Id - \Pi\cdot\mS$, the operator
$\mS$ is block-diagonal, each block being associated to a different
subdomain, so that matrix-vector product is trivially parellelisable.
Of course, each block of $\mS$ involves a DtN operator.

The main new feature of the formulation we present here is the transmission
operator $\Pi$. Contrary to the exchange operator traditionally used in
OSM, see e.g. Formula (42) in \cite{MR1764190}, our transmission
operator $\Pi$ is not local anymore. But it only involves exponentially
decaying kernels, with a damping factor $\gamma$ that can be tuned, so 
that $\Pi$ can nevertheless be considered quasi-local. In addition,
various techniques (H-matrices \cite{MR2451321,MR2767920,MR3445676},
Fast Multipole Method \cite{MR1756765,MR1489257}) can be used to sparsify
this operator further.

\subsection{Well-posedness of the new formulation}
Let us examine the properties of the operator $\Id-\Pi\cdot\mS$ in
detail. First of all $\Pi\cdot\mS$ continuously maps $\mbH_{\neu}(\Gamma)$
into $\mbH_{\neu}(\Gamma)$. In addition, combining Corollary \ref{ContractionProj}
and Proposition \ref{Contractivity}, we obtain a contractivity result.
\begin{lemma}\label{GlobalContractivity}\quad\\
  We have $\Vert \Pi\cdot\mS(\ctrp)\Vert_{\mbH_{\neu}(\Gamma)}\leq
  \Vert\ctrp\Vert_{\mbH_{\neu}(\Gamma)}$ for all $\ctrp\in\mbH_{\neu}(\Gamma)$.
\end{lemma}

\noindent 
A direct consequence of this property is that the numerical range
of the operator $\Id - \Pi\cdot\mS$ is located in the complex right-half
plane $\mathbb{C}_{+}:=\{z\in\mathbb{C}, \;\Re e\{z\}\geq 0\}$. This is
definitely an interesting feature from the perspective of linear solvers.
Next this operator is also one-to-one.

\begin{proposition}\label{InjectivityProp}\quad\\
  $\mrm{ker}(\Id - \Pi\cdot\mS) = \{0\}$.
\end{proposition}
\noindent \textbf{Proof:}

Consider a $\ctrp = (p^{j})_{j=0}^{\mJ}\in\mbH_{\neu}(\Gamma)$ satisfying $\ctrp = \Pi\mS(\ctrp)$.
Consider the function $v\in\mL^{2}_{\loc}(\RR^{d})$ such that, its restriction in each subdomain
$v\vert_{\Omega_{j}}$ belongs to $\mH^{1}_{\loc}(\overline{\Omega}_{j})$ and satisfies
$\div(\mu\nabla v)+\kappa^{2}v = 0$ in $\Omega_{j}$, $v$ is $\kappa_{0}$-outgoing and
$\tau_{-}^{j}(v) = p^{j}$. By construction we have $\tau_{-}(v) = \ctrp$ and $\tau_{+}(v) = \mS(\ctrp)$.
Setting $\ctrv = (\ctrv_{\dir},\ctrv_{\neu}):=(\tau^{j}_{\dir}(v),\mu_{j}
\tau^{j}_{\neu}(v))_{j=0\dots \mJ}$, we have $0 = \ctrp - \Pi\mS(\ctrp) =
\ctrv_{\neu}-\imath\omega \mT(\ctrv_{\dir}) -
\Pi(\ctrv_{\neu}+\imath\omega\mT(\ctrv_{\dir}))$. Hence, applying 
Proposition \ref{CaracTransmissionRobinTraces}, we deduce that $\ctrv\in\mbX(\Gamma)$. 
Since, on the other hand, we have  $\ctrv\in\mathscr{C}^{+}(\Gamma)$ by construction, 
we conclude that  $\ctrv\in\mathscr{C}^{+}(\Gamma)\cap \mbX(\Gamma) = \{0\}$
according to Proposition \ref{DecompSTFPropagCauchy}.
Hence $\ctrp = \ctrv_{\neu}-\imath\omega\mT(\ctrv_{\dir}) = 0$.
\hfill $\Box$

\quad\\
The operator $\Id - \Pi\mS$ is actually coercive. 

\begin{theorem}\label{CoercivityThm}\quad\\
  There exists $\alpha>0$ such that
  $\Re e\{((\Id - \Pi\cdot\mS)\ctrp,\ctrp)_{\mbH_{\neu}(\Gamma)}\}\geq \alpha
  \Vert \ctrp\Vert_{\mbH_{\neu}(\Gamma)}^{2}$ for all $\ctrp\in \mbH_{\neu}(\Gamma)$. 
\end{theorem}
\noindent \textbf{Proof:}

We need first to introduce a few notations that we shall use only for this proof.
According to Proposition \ref{DecompSTFPropagCauchy}, there exists
a bounded projection operator $\mQ:\mbH(\Gamma)\to \mbH(\Gamma)$
with $\mrm{range}(\mQ) = \mathscr{C}^{+}(\Gamma)$ and $\mrm{ker}(\mQ) = \mbX(\Gamma)$.
For convenience, we set
\begin{equation}\label{defscalednorm}
  \begin{array}{l}
    \Vert \mQ\Vert_{\omega} :=  \sup_{\ctrv\in \mbH(\Gamma)\setminus\{0\} }
    \Vert \mQ(\ctrv)\Vert_{\omega}/
    \Vert \ctrv\Vert_{\omega}\\[5pt]
    \text{where}\;\;
    \Vert \ctrv\Vert_{\omega}^{2} :=
    \Vert \ctrv_{\neu}\Vert_{\mbH_{\neu}(\Gamma)}^{2} +
    \omega^{2}\Vert \ctrv_{\dir}\Vert_{\mbH_{\dir}(\Gamma)}^{2}
  \end{array}
\end{equation}
Because $\omega>0$ is a simple fixed positive constant,
$\Vert \;\Vert_{\omega}$ and $\Vert \;\Vert_{\mbH(\Gamma)}$ are equivalent norms,
and continuity of the projection $\mQ$ is exactly equivalent to the boundedness
of $\Vert \mQ\Vert_{\omega}$. We shall also consider the bounded orthogonal
projectors $\mP_{\pm}:\mbH_{\neu}(\Gamma)\to \mbH_{\neu}(\Gamma)$ defined by
\begin{equation}
  \mP_{\pm} = (\Id\pm\Pi)/2
\end{equation}
Now pick an arbitrary $\ctrp\in \mbH_{\neu}(\Gamma)$. 
Set $\ctrf := (\Id - \Pi\mS)\ctrp$, and 
define $\ctrg_{\dir}:=\imath\omega^{-1}\mT^{-1}(\Id+\Pi)\ctrf/4$ and
$\ctrg_{\neu}:=(\Id-\Pi)\ctrf/4$ and $\ctrg:=(\ctrg_{\dir},\ctrg_{\neu})\in\mbH(\Gamma)$.
The tuple of traces $\ctru = \mQ(\ctrg)\in\mathscr{C}^{+}(\Gamma)$ satisfies
$\ctrg-\ctru\in\mbX(\Gamma)$ so, applying Proposition \ref{CaracTransmissionRobinTraces},
we also have $\ctru_{\neu}-\ctrg_{\neu}-\imath\omega\mT(\ctru_{\dir}-\ctrg_{\dir}) =
\Pi(\ctru_{\neu}-\ctrg_{\neu}+\imath\omega\mT(\ctru_{\dir}-\ctrg_{\dir}) )$
which rewrites
\begin{equation}\label{ConstructInterm2}
  \begin{array}{l}
    \ctru_{\neu}-\imath\omega\mT(\ctru_{\dir}) -
    \Pi(\ctru_{\neu}+\imath\omega\mT(\ctru_{\dir}))\\[5pt]
    \hspace{1cm}=(\Id - \Pi)\ctrg_{\neu} - \imath\omega(\Id + \Pi)\mT(\ctrg_{\dir})\\[5pt]
    \hspace{1cm}=\mP_{-}^{2}\ctrf + \mP_{+}^{2}\ctrf = (\mP_{-}+ \mP_{+})\ctrf = \ctrf
  \end{array}
\end{equation}
Due to the continuity of $\mQ$, we obviously have
$\Vert \ctru\Vert_{\omega}\leq \Vert \mQ\Vert_{\omega}\cdot
\Vert \ctrg\Vert_{\omega}$, where $\Vert \mQ\Vert_{\omega}$ is defined
with (\ref{defscalednorm}).  On the other hand multiplying (\ref{ConstructInterm2})
on the left by $\mP_{\pm}$ we obtain 
\begin{equation}\label{Isometry2}
  \begin{aligned}
    & \mP_{+}\mT(\ctru_{\dir}) = \imath\omega^{-1}\mP_{+}(\ctrf)/2 = \mT(\ctrg_{\dir})\\
    & \mP_{-}(\ctru_{\neu}) = \mP_{-}(\ctrf)/2 = \ctrg_{\neu}\\
    & \Rightarrow\quad \Vert \ctrg\Vert_{\omega}^{2}=
    \omega^{2}\Vert \mP_{+}\mT(\ctru_{\dir})\Vert_{\mbH_{\neu}(\Gamma)}^{2} +
    \Vert \mP_{-}(\ctru_{\neu})\Vert_{\mbH_{\neu}(\Gamma)}^{2}\\
  \end{aligned}
\end{equation}
which shows that $\Vert \ctru\Vert_{\omega}\leq \Vert \mQ\Vert_{\omega}
(\omega^{2}\Vert \mP_{+}\mT(\ctru_{\dir})\Vert_{\mbH_{\neu}(\Gamma)}^{2} +
\Vert \mP_{-}(\ctru_{\neu})\Vert_{\mbH_{\neu}(\Gamma)}^{2})$. Next
observe that (\ref{ConstructInterm2}) implies
$(\Id - \Pi\mS)(\ctru_{\neu}-\imath\omega\mT(\ctru_{\dir}))
= \ctrf$ hence, according to Proposition \ref{InjectivityProp}, $\ctrp =
\ctru_{\neu}-\imath\omega\mT(\ctru_{\dir})$, which leads to the estimate
\begin{equation}\label{EstimateCoercivity1}
  \Vert \ctrp\Vert^{2}_{\mbH_{\neu}(\Gamma)}/2
  \leq \Vert \ctru\Vert_{\omega}^{2} \leq
  \Vert \mQ\Vert_{\omega}^{2} \Vert\ctrg\Vert_{\omega}^{2}.
\end{equation}
Since the projectors $\mP_{\pm}$ are orthogonal for the scalar
 product $(\cdot,\cdot)_{\mbH_{\neu}(\Gamma)}$ we obtain 
\begin{equation}
  \begin{aligned}
    & \frac{1}{2}(\ctrp-\Pi\mS(\ctrp),\ctrp)_{\mbH_{\neu}(\Gamma)}\\
    & \quad = \frac{1}{2}(\ctru_{\neu}-\imath \omega\mT(\ctru_{\dir})-
    \Pi(\ctru_{\neu}+\imath\omega \mT(\ctru_{\dir})),
    \ctru_{\neu}-\imath \omega\mT(\ctru_{\dir}))_{\mbH_{\neu}(\Gamma)}\\
    & \quad = (\mP_{-}(\ctru_{\neu})-
    \imath\omega\mP_{+}\mT(\ctru_{\dir}),\ctru_{\neu}-\imath
    \omega\mT(\ctru_{\dir}))_{\mbH_{\neu}(\Gamma)}\\
    & \quad = \Vert \mP_{-}(\ctru_{\neu})\Vert_{\mbH_{\neu}(\Gamma)}^{2}
    + \omega^{2}\Vert \mP_{+}\mT(\ctru_{\dir})\Vert_{\mbH_{\neu}(\Gamma)}^{2}\\
    & \quad - \imath \omega(\mP_{+}\mT(\ctru_{\dir}),\ctru_{\neu} )_{\mbH_{\neu}(\Gamma)}
    + \imath \omega(\mP_{-}(\ctru_{\neu}),\mT(\ctru_{\dir}) )_{\mbH_{\neu}(\Gamma)}
  \end{aligned}
\end{equation}
Using the identity obtained in (\ref{Isometry2}) to replace
$\Vert \mP_{-}(\ctru_{\neu})\Vert_{\mbH_{\neu}(\Gamma)}^{2}
+ \omega^{2}\Vert \mP_{+}\mT(\ctru_{\dir})\Vert_{\mbH_{\neu}(\Gamma)}^{2}$
in the identity above, using that $\mP_{\pm} = (\Id\pm\Pi)/2$, and observing
that $(\mT(\ctrv),\ctrp)_{\mbH_{\neu}(\Gamma)} = \llangle \ctrv,\overline{\ctrp}\rrangle$,
we obtain
\begin{equation}
  \begin{aligned}
    & \frac{1}{2}(\ctrp-\Pi\mS(\ctrp),\ctrp)_{\mbH_{\neu}(\Gamma)}\\
    & \quad = \Vert \ctrg\Vert_{\omega}^{2}-
    (\imath/2) \llangle\ctru_{\dir},\overline{\ctru}_{\neu}\rrangle +
    (\imath/2)\llangle\ctru_{\neu},\overline{\ctru}_{\dir}\rrangle\\
    & \quad -
    (\imath/2)(\Pi\mT(\ctru_{\dir}),\ctru_{\neu})_{\mbH_{\neu}(\Gamma)} -
    (\imath/2)(\Pi(\ctru_{\neu}),\mT(\ctru_{\dir}))_{\mbH_{\neu}(\Gamma)}\\
    & \quad =  \Vert \ctrg\Vert_{\omega}^{2} -(\imath/2)\llbr \ctru,\overline{\ctru}\rrbr
    -\imath \Re e\{ (\Pi\mT(\ctru_{\dir}),\ctru_{\neu})_{\mbH_{\neu}(\Gamma)}\}
  \end{aligned}
\end{equation}
Using Lemma \ref{EnergyFlux},
the real part of the previous identity is bounded from below by 
$\Re e\{(\ctrp-\Pi\mS(\ctrp),\ctrp)_{\mbH_{\neu}(\Gamma)}\}\geq
2\Vert \ctrg\Vert_{\omega}^{2}$. We conclude by using
(\ref{EstimateCoercivity1}).  \hfill $\Box$

\quad\\
Lax-Milgram lemma combined with the previous theorem yields bijectivity
of $\Id - \Pi\mS$ as an obvious outcome.

\begin{corollary}\quad\\
  The operator $\Id - \Pi\mS:\mbH_{\neu}(\Gamma)\to \mbH_{\neu}(\Gamma)$
  is an isomorphism.
\end{corollary}

\subsection{Solution strategy}
Let us briefly discuss how, in practice, to solve (\ref{OSMFormulation}) i.e.
an equation of the form $\ctrp - \Pi\cdot\mS(\ctrp) = \ctrf$. First of all,
since $\Pi^{2} = \Id$, this equation can be transformed into 
$(\Pi - \mS)\ctrp = \Pi(\ctrf) = \tau_{+}(\phi_{f})$
which is practically  more convenient as it avoids handling a product
of operators. A general Krylov solver such as GMRes could be considered
for solving this equation. We refer the reader to \cite[chap.6]{zbMATH01953444}
for more details on this solver.

\paragraph{Convergence of Richardson's linear solver}
An alternative more straightforward strategy relies on Richardson's
iterative method \cite[chap.6]{zbMATH01953444}, \cite[\S 9.1]{zbMATH05189473}
that writes 
\begin{equation}
  \ctrp^{n+1} =  (1-\beta)\ctrp^{(n)} +\beta\Pi\mS\cdot\ctrp^{(n)} + \beta\ctrf
\end{equation}
where $\beta\in (0,1)$ is a relaxation parameter. Following Theorem 7 and
Remark 9 in \cite{MR1764190}, a rough estimate can be derived for the convergence
of Richardson's linear solver in this case. Let $\ctrp^{\infty}$
refer to the unique solution to (\ref{OSMFormulation}) and set
$\mathfrak{e}^{(n)}:= \ctrp^{\infty} - \ctrp^{(n)}$ so that
$\mathfrak{e}^{(n+1)} = ((1-\beta)\Id + \beta\Pi\mS) \mathfrak{e}^{(n)}$.
Recall the convexity identity
\begin{equation}\label{ConvexityIdentity}
  \begin{aligned}
    \Vert (1-\beta)\mathbf{x} + \beta \mathbf{y}
    \Vert_{\mbH_{\neu}(\Gamma)}^{2}
    & = (1-\beta)\Vert\mathbf{x}\Vert_{\mbH_{\neu}(\Gamma)}^{2} +
    \beta \Vert\mathbf{y}\Vert_{\mbH_{\neu}(\Gamma)}^{2}\\
    & \textcolor{white}{=}
    -\beta(1-\beta) \Vert\mathbf{x-y}\Vert_{\mbH_{\neu}(\Gamma)}^{2}
  \end{aligned}
\end{equation}
which holds for any $\mathbf{x},\mathbf{y}\in \mbH_{\neu}(\Gamma)$
and any $\beta\in (0,1)$. In addition the coercivity
estimate of Theorem \ref{CoercivityThm} yields the lower bound 
$\Vert (\Id - \Pi\mS)\ctrp\Vert_{\mbH_{\neu}(\Gamma)}\geq
\alpha \Vert\ctrp\Vert_{\mbH_{\neu}(\Gamma)}\forall\ctrp\in \mbH_{\neu}(\Gamma)$.
Combining this lower bound with Lemma \ref{GlobalContractivity} and
(\ref{ConvexityIdentity}) thus yields
\begin{equation*}
  \begin{aligned}
    \Vert \mathfrak{e}^{(n+1)}\Vert_{\mbH_{\neu}(\Gamma)}^{2}
    & = \Vert (1-\beta)\mathfrak{e}^{(n)} +\beta
    \Pi\mS\cdot \mathfrak{e}^{(n)}\Vert_{\mbH_{\neu}(\Gamma)}^{2}\\
    & =  (1-\beta)\Vert\mathfrak{e}^{(n)}\Vert_{\mbH_{\neu}(\Gamma)}^{2}
    +\beta \Vert \Pi\mS\cdot \mathfrak{e}^{(n)}\Vert_{\mbH_{\neu}(\Gamma)}^{2}\\
    & \textcolor{white}{=}
    - \beta(1-\beta)\Vert (\Id - \Pi\mS)\mathfrak{e}^{(n)}\Vert_{\mbH_{\neu}(\Gamma)}^{2}\\
    & \leq (1-\alpha^{2}\beta(1-\beta))
    \Vert\mathfrak{e}^{(n)}\Vert_{\mbH_{\neu}(\Gamma)}^{2}
  \end{aligned}
\end{equation*}
In this estimate, the convergence factor $(1-\alpha^{2}\beta(1-\beta))^{1/2}<1$
is thus minimized for $\beta = 1/2$ and takes the value $(1-(\alpha/2)^{2})^{1/2}$
in this case.

\section*{Acknowledgement}
This work received support from the French National Research Agency (ANR)
through the NonlocalDD project, grant ref. ANR-15-CE23-0017-01.


\begin{thebibliography}{10}

\bibitem{zbMATH05189473}
G.~Allaire and S.M. Kaber.
\newblock {\em {Numerical linear algebra.}}, volume~55.
\newblock New York, Springer, 2008.

\bibitem{zbMATH06732120}
X.~{Antoine} and C.~{Geuzaine}.
\newblock {Optimized Schwarz domain decomposition methods for scalar and vector
  Helmholtz equations.}
\newblock In {\em {Modern solvers for Helmholtz problems}}, pages 189--213.
  Basel: Birkh\"auser/Springer, 2017.

\bibitem{MR2451321}
M.~Bebendorf.
\newblock {\em Hierarchical matrices}, volume~63 of {\em Lecture Notes in
  Computational Science and Engineering}.
\newblock Springer-Verlag, Berlin, 2008.
\newblock A means to efficiently solve elliptic boundary value problems.

\bibitem{MR2767920}
S.~B\"{o}rm.
\newblock {\em Efficient numerical methods for non-local operators}, volume~14
  of {\em EMS Tracts in Mathematics}.
\newblock European Mathematical Society (EMS), Z\"{u}rich, 2010.
\newblock ${\mathscr{H}}^2$-matrix compression, algorithms and analysis.

\bibitem{zbMATH06039354}
Y.~{Boubendir}, X.~{Antoine}, and C.~{Geuzaine}.
\newblock {A quasi-optimal non-overlapping domain decomposition algorithm for
  the Helmholtz equation.}
\newblock {\em {J. Comput. Phys.}}, 231(2):262--280, 2012.

\bibitem{C11_144}
X.~Claeys.
\newblock A single trace integral formulation of the second kind for acoustic
  scattering.
\newblock Technical Report 2011-14, Seminar for Applied Mathematics, ETH
  Z{\"u}rich, Switzerland, 2011.

\bibitem{MR3403719}
X.~Claeys.
\newblock Quasi-local multitrace boundary integral formulations.
\newblock {\em Numer. Methods Partial Differential Equations},
  31(6):2043--2062, 2015.

\bibitem{zbMATH06185811}
X.~Claeys and R.~Hiptmair.
\newblock {Electromagnetic scattering at composite objects: a novel multi-trace
  boundary integral formulation.}
\newblock {\em {ESAIM, Math. Model. Numer. Anal.}}, 46(6):1421--1445, 2012.

\bibitem{MR3101780}
X.~Claeys and R.~Hiptmair.
\newblock Integral equations on multi-screens.
\newblock {\em Integral Equations Operator Theory}, 77(2):167--197, 2013.

\bibitem{MR3069956}
X.~Claeys and R.~Hiptmair.
\newblock Multi-trace boundary integral formulation for acoustic scattering by
  composite structures.
\newblock {\em Comm. Pure Appl. Math.}, 66(8):1163--1201, 2013.

\bibitem{zbMATH06286985}
X.~Claeys, R.~Hiptmair, and C.~Jerez-Hanckes.
\newblock {Multitrace boundary integral equations.}
\newblock In {\em {Direct and inverse problems in wave propagation and
  applications. Selected papers of the workshop on wave propagation and
  scattering, inverse problems and applications in energy and the environment,
  Linz, Austria, November 21--25, 2011}}, pages 51--100. Berlin: de Gruyter,
  2013.

\bibitem{MR3313601}
X.~Claeys, R.~Hiptmair, and E.~Spindler.
\newblock A second-kind {G}alerkin boundary element method for scattering at
  composite objects.
\newblock {\em BIT}, 55(1):33--57, 2015.

\bibitem{MR3720391}
X.~Claeys, R.~Hiptmair, and E.~Spindler.
\newblock Second kind boundary integral equation for multi-subdomain diffusion
  problems.
\newblock {\em Adv. Comput. Math.}, 43(5):1075--1101, 2017.

\bibitem{MR3725822}
X.~Claeys, R.~Hiptmair, and E.~Spindler.
\newblock Second-kind boundary integral equations for electromagnetic
  scattering at composite objects.
\newblock {\em Comput. Math. Appl.}, 74(11):2650--2670, 2017.

\bibitem{10.1007/978-3-319-93873-8_16}
X.~Claeys, B.~Thierry, and F.~Collino.
\newblock Integral equation based optimized schwarz method for
  electromagnetics.
\newblock In Petter~E. Bj{\o}rstad, Susanne~C. Brenner, Lawrence Halpern,
  Hyea~Hyun Kim, Ralf Kornhuber, Talal Rahman, and Olof~B. Widlund, editors,
  {\em Domain Decomposition Methods in Science and Engineering XXIV}, pages
  187--194, Cham, 2018. Springer International Publishing.

\bibitem{MR1764190}
F.~Collino, S.~Ghanemi, and P.~Joly.
\newblock Domain decomposition method for harmonic wave propagation: a general
  presentation.
\newblock {\em Comput. Methods Appl. Mech. Engrg.}, 184(2-4):171--211, 2000.

\bibitem{MR2986407}
D.~Colton and R.~Kress.
\newblock {\em Inverse acoustic and electromagnetic scattering theory},
  volume~93 of {\em Applied Mathematical Sciences}.
\newblock Springer, New York, third edition, 2013.

\bibitem{MR1756765}
E.~Darve.
\newblock The fast multipole method: numerical implementation.
\newblock {\em J. Comput. Phys.}, 160(1):195--240, 2000.

\bibitem{MR1071633}
B.~Despr\'{e}s.
\newblock D\'{e}composition de domaine et probl\`eme de {H}elmholtz.
\newblock {\em C. R. Acad. Sci. Paris S\'{e}r. I Math.}, 311(6):313--316, 1990.

\bibitem{MR1105979}
B.~Despr\'{e}s.
\newblock Domain decomposition method and the {H}elmholtz problem.
\newblock In {\em Mathematical and numerical aspects of wave propagation
  phenomena ({S}trasbourg, 1991)}, pages 44--52. SIAM, Philadelphia, PA, 1991.

\bibitem{MR1291197}
B.~Despr\'{e}s.
\newblock {\em M\'{e}thodes de d\'{e}composition de domaine pour les
  probl\`emes de propagation d'ondes en r\'{e}gime harmonique. {L}e
  th\'{e}or\`eme de {B}org pour l'\'{e}quation de {H}ill vectorielle}.
\newblock Institut National de Recherche en Informatique et en Automatique
  (INRIA), Rocquencourt, 1991.
\newblock Th\`ese, Universit\'{e} de Paris IX (Dauphine), Paris, 1991.

\bibitem{MR1227838}
B.~Despr\'{e}s.
\newblock Domain decomposition method and the {H}elmholtz problem. {II}.
\newblock In {\em Second {I}nternational {C}onference on {M}athematical and
  {N}umerical {A}spects of {W}ave {P}ropagation ({N}ewark, {DE}, 1993)}, pages
  197--206. SIAM, Philadelphia, PA, 1993.

\bibitem{zbMATH06534518}
V.~Dolean, P.~Jolivet, and F.~Nataf.
\newblock {\em {An introduction to domain decomposition methods. Algorithms,
  theory, and parallel implementation.}}
\newblock Philadelphia, PA: Society for Industrial and Applied Mathematics
  (SIAM), 2015.

\bibitem{zbMATH06666466}
M.~{El Bouajaji}, X.~{Antoine}, and C.~{Geuzaine}.
\newblock {Approximate local magnetic-to-electric surface operators for
  time-harmonic Maxwell's equations.}
\newblock {\em {J. Comput. Phys.}}, 279:241--260, 2014.

\bibitem{zbMATH06660584}
M.~{El Bouajaji}, B.~{Thierry}, X.~{Antoine}, and C.~{Geuzaine}.
\newblock {A quasi-optimal domain decomposition algorithm for the time-harmonic
  Maxwell's equations.}
\newblock {\em {J. Comput. Phys.}}, 294:38--57, 2015.

\bibitem{zbMATH07020343}
M.J. Gander and H.~Zhang.
\newblock {A class of iterative solvers for the Helmholtz equation:
  factorizations, sweeping preconditioners, source transfer, single layer
  potentials, polarized traces, and optimized Schwarz methods.}
\newblock {\em {SIAM Rev.}}, 61(1):3--76, 2019.

\bibitem{MR1489257}
L.~Greengard and V.~Rokhlin.
\newblock A new version of the fast multipole method for the {L}aplace equation
  in three dimensions.
\newblock In {\em Acta numerica, 1997}, volume~6 of {\em Acta Numer.}, pages
  229--269. Cambridge Univ. Press, Cambridge, 1997.

\bibitem{MR3445676}
W.~Hackbusch.
\newblock {\em Hierarchical matrices: algorithms and analysis}, volume~49 of
  {\em Springer Series in Computational Mathematics}.
\newblock Springer, Heidelberg, 2015.

\bibitem{zbMATH02038442}
U.~{Langer} and O.~{Steinbach}.
\newblock {Boundary element tearing and interconnecting methods.}
\newblock {\em {Computing}}, 71(3):205--228, 2003.

\bibitem{LecouvezThesis}
M.~Lecouvez.
\newblock {\em {Iterative methods for domain decomposition without overlap with
  exponential convergence for the Helmholtz equation}}.
\newblock Phd thesis, {Ecole Polytechnique}, July 2015.

\bibitem{LECOUVEZ2014403}
M.~Lecouvez, B.~Stupfel, P.~Joly, and F.~Collino.
\newblock Quasi-local transmission conditions for non-overlapping domain
  decomposition methods for the helmholtz equation.
\newblock {\em Comptes Rendus Physique}, 15(5):403 -- 414, 2014.
\newblock Electromagnetism / Électromagnétisme.

\bibitem{MR841971}
R.~Leis.
\newblock {\em Initial-boundary value problems in mathematical physics}.
\newblock B. G. Teubner, Stuttgart; John Wiley \& Sons, Ltd., Chichester, 1986.

\bibitem{modave:hal-01925160}
Axel Modave, Christophe Geuzaine, and Xavier Antoine.
\newblock {Corner treatment for high-order local absorbing boundary conditions
  in high-frequency acoustic scattering}.
\newblock working paper or preprint, November 2018.

\bibitem{zbMATH06308881}
A.~Moiola and E.A. Spence.
\newblock {Is the Helmholtz equation really sign-indefinite?}
\newblock {\em {SIAM Rev.}}, 56(2):274--312, 2014.

\bibitem{MR2723248}
F.W.J. Olver, D.W. Lozier, R.~F. Boisvert, and C.W. Clark, editors.
\newblock {\em N{IST} handbook of mathematical functions}.
\newblock U.S. Department of Commerce, National Institute of Standards and
  Technology, Washington, DC; Cambridge University Press, Cambridge, 2010.

\bibitem{zbMATH06029154}
C.~Pechstein.
\newblock {\em {Finite and boundary element tearing and interconnecting solvers
  for multiscale problems.}}, volume~90.
\newblock Berlin: Springer, 2013.

\bibitem{zbMATH01953444}
Y.~Saad.
\newblock {\em {Iterative methods for sparse linear systems. 2nd ed.}}
\newblock Philadelphia, PA: SIAM Society for Industrial and Applied
  Mathematics, 2nd ed. edition, 2003.

\bibitem{zbMATH02113718}
A.~Toselli and O.~Widlund.
\newblock {\em {Domain decomposition methods -- algorithms and theory.}},
  volume~34.
\newblock Berlin: Springer, 2005.

\bibitem{zbMATH07007340}
A.~{Vion} and C.~{Geuzaine}.
\newblock {Improved sweeping preconditioners for domain decomposition
  algorithms applied to time-harmonic Helmholtz and Maxwell problems.}
\newblock {\em {ESAIM, Proc. Surv.}}, 61:93--111, 2018.

\bibitem{zbMATH04136456}
T.~{von Petersdorff}.
\newblock {Boundary integral equations for mixed Dirichlet, Neumann and
  transmission problems.}
\newblock {\em {Math. Methods Appl. Sci.}}, 11(2):185--213, 1989.

\end{thebibliography}

\end{document}